\documentclass[10pt, reqno]{amsart}

\usepackage{amssymb,latexsym,amsmath,amsfonts}
\usepackage{amsmath,amscd,amsthm,amssymb}

\numberwithin{equation}{section}

\newcommand{\Zz}{{\mathbb Z}}
\newcommand{\Ff}{{\mathbb F}}
\newcommand{\Cc}{{\mathbb C}}
\newcommand{\Qq}{{\mathbb Q}}

\newcommand{\p}{{\mathbf P}}
\newcommand{\q}{{\mathbf Q}}

\newcommand{\M}{{\mathbf M}}
\newcommand{\A}{{\mathbf A}}
\newcommand{\B}{{\mathbf B}}
\newcommand{\C}{{\mathbf C}}
\newcommand{\D}{{\mathbf D}}
\newcommand{\h}{{\mathcal H}}
\newcommand{\x}{{\mathbf x}}

\newcommand{\Proof}{ \noindent{\bf Proof:}\quad }

\def\QED{\qed\medskip\par}
\newtheorem{Theorem} {Theorem} [section]
\newtheorem{Proposition} [Theorem] {Proposition}
\newtheorem{Lemma} [Theorem] {Lemma}
\newtheorem{Corollary} [Theorem] {Corollary}
\newtheorem{Conjecture}[Theorem]{Conjecture}

\newtheorem{Definition}[Theorem]{Definition}
\newtheorem{Remark} [Theorem] {Remark}

\def\Sq{\;\Box}
\def\Nsq{\;\not\!\!\Box}

\begin{document}

\title[Dimensions of Binary Codes]
{Dimensions of Some Binary Codes Arising From A Conic in $PG(2,q)$}

\author[Sin, Wu and Xiang]{Peter Sin, Junhua Wu and Qing Xiang$^{\dagger}$}

\thanks{$^{\dagger}$Research supported in part by NSF Grant DMS 0701049, and by the
Overseas Cooperation Fund (grant 10928101) of China.}
\address{Department of Mathematics, University of Florida, Gainesville, FL 32611, USA}
\email{sin@ufl.edu}

\address{Department of Mathematical Sciences, Worcester Polytechnic Institute, Worcester, MA 01609, USA}
\email{wuj@wpi.edu}

\address{Department of Mathematical Sciences, University of Delaware, Newark, DE 19716, USA} \email{xiang@math.udel.edu}


\begin{abstract}
Let $\mathcal{O}$ be a conic in the classical projective plane
$PG(2,q)$, where $q$ is an odd prime power. With respect to
$\mathcal{O}$, the lines of $PG(2,q)$ are classified as passant,
tangent, and secant lines, and the points of $PG(2,q)$ are
classified as internal, absolute and external points. The incidence
matrices between the secant/passant lines and the external/internal
points were used in \cite{keith1} to produce several classes of
structured low-density parity-check binary codes. In particular, the
authors of \cite{keith1} gave conjectured dimension formula for the
binary code $\mathcal{L}$ which arises as the $\Ff_2$-null space of
the incidence matrix between the secant lines and the external
points to $\mathcal{O}$. In this paper, we prove the conjecture on
the dimension of $\mathcal{L}$ by using a combination of techniques
from finite geometry and modular representation theory.\end{abstract}

\maketitle

\section{Introduction}\label{intro}

Let $\Ff_q$ be the finite field of order $q$, where $q=p^e$, $p$ is
a prime and $e\geq 1$ is an integer. Let $PG(2,q)$ denote the
classical projective plane of order $q$ constructed from the
3-dimensional vector space $\Ff_q^3$ in the standard way. A {\it
conic} in $PG(2,q)$ is the set of points $(x,y,z)$
satisfying a non-zero quadratic form. We say that a conic is {\it
non-degenerate} if it does not contain an entire line of $PG(2,q)$.
By a linear change of coordinates, any non-degenerate conic is
equivalent to
\begin{equation}\label{oval}
\mathcal{O} = \{(1,t,t^2)\mid t\in \Ff_q\}\cup\{(0,0,1)\},
\end{equation}
the set of (projective) $\Ff_q$-zeros of the non-degenerate
quadratic form
\begin{equation}\label{NQ}
Q(X_0,X_1,X_2) = X_1^2-X_0 X_2
\end{equation}
over $\Ff_q$.

It can be shown \cite[p.~157]{hir} that every non-degenerate conic
has $q+1$ points, no three of which are collinear. That is, a
non-degenerate conic is an oval. When $q$ is odd, Segre \cite{segre}
proved that an oval in $PG(2,q)$ must be a non-degenerate conic. It
follows that in $PG(2,q)$, where $q$ is odd, ovals and
non-degenerate conics are the same objects.

In the rest of this paper, we will always assume that $q=p^e$ is an
{\it odd} prime power, and fix the conic in (\ref{oval}) as the
``standard" conic. A line $\ell$ is called a {\it passant}, a {\it
tangent}, or a {\it secant} of $\mathcal{O}$ according as
$|\ell\cap\mathcal{O}|=0$, $1$, or $2$. Since the conic
$\mathcal{O}$ is an oval, we see that every line of $PG(2,q)$ falls
into one of these classes. A point $\p$ is called an {\it internal},
{\it absolute}, or {\it external} point according as $\p$ lies on
$0$, $1$, or $2$ tangent lines to $\mathcal{O}$. It is an easy
exercise to show that in $PG(2,q)$, $\mathcal{O}$ has
$\frac{q(q+1)}{2}$ secant lines, $q+1$ tangent lines, and
$\frac{q(q-1)}{2}$ passant lines; $\frac{q(q+1)}{2}$ external
points, $q+1$ absolute points (which are the points on ${\mathcal
O}$) and $\frac{q(q-1)}{2}$ internal points. We will denote the sets
of secant, tangent, and passant lines by $Se$, $T$ and $Pa$,
respectively, and the sets of external and internal points by $E$
and $I$, respectively. In fact, the quadratic form $Q$ in (\ref{NQ})
induces a polarity $\perp$ (a correlation of order 2) of $PG(2,q)$
under which $E$ and $Se$, ${\mathcal O}$ and $T$, and $I$ and $Pa$
are in one-to-one correspondence with each other, respectively. A
summary of the intersection patterns for the various types of points
and lines is given in Tables~\ref{tab1} and \ref{tab2} in Section 2.

Let $\A$ be the $(q^2+q+1)$ times $(q^2+q+1)$ line-point incidence
matrix of $PG(2,q)$. That is, the rows and columns of $\A$ are
labeled by the lines and points of $PG(2,q)$, respectively, and the
$(\ell,\p)$-entry of $\A$ is $1$ if $\p\in\ell$, $0$ otherwise. It
is well known that the $2$-rank of $\A$ is $q^2+q$ ({\cite{HP}}) and
the $p$-rank of $\A$ is $\binom{p+1}{2}^e+1$ ({\cite{AK}}), where
$q=p^e$.

In \cite{keith1}, Droms, Mellinger and Meyer considered the
following partition of $\A$ into nine submatrices:
\begin{equation}\label{matrix_A}
\A=
\left(\begin{array}{lll}
\A_{11} & \A_{12} & \A_{13} \\
\A_{21} & \A_{22} & \A_{23} \\
\A_{31} & \A_{32} & \A_{33}\\
\end{array}
\right),
\end{equation}
where the rows of $\A_{11}$, $\A_{21}$, and $\A_{31}$ are labeled by
the tangent, passant, and secant lines respectively, and the columns
of $\A_{11}$, $\A_{12}$, and $\A_{13}$ are labeled by the absolute,
internal, and external points, respectively. These authors used the
submatrices $\A_{ij}$ for $2\le i,j\le 3$ to construct four binary
linear codes, and showed that these codes are good examples of
structured low-density parity-check (LDPC) codes. Based on
computational evidence, they made conjectures on the dimensions of
these binary LDPC codes. In particular, it was conjectured in
\cite{keith1} that the dimension of the $\Ff_2$-null space of
$\A_{33}$ (i.e. the incidence matrix of secant lines versus external
points) is given by the following simple formula.
\begin{Conjecture}\label{con1}
Let $\mathcal{L}$ be the $\Ff_2$-null space of $\A_{33}$. Then
\begin{displaymath}
{\rm dim}_{\Ff_2}(\mathcal{L})=
\begin{cases}
\frac{(q-1)^2}{4}+1, &\text{if}\;q\equiv 1\pmod 4,\\
\frac{(q-1)^2}{4}-1, &\text{if}\;q\equiv 3\pmod 4.\\
\end{cases}
\end{displaymath}
\end{Conjecture}

In this paper, we use a combination of techniques from finite
geometry and group representation theory to confirm the above
conjecture. The conjectured formulas for the dimensions of the
binary codes arising from $\A_{22}$, $\A_{23}$ and $\A_{32}$ will be
proved in forthcoming papers. We remark in passing that the
$p$-ranks of all submatrices in (\ref{matrix_A}) were recently
computed in {\cite{wu}}, where $p$ is the characteristic of the
defining field of $PG(2,q)$ and $q=p^e$. For instance,
$rank_p(\A_{33})=\binom{p+1}{2}^e$ {\cite[Theorem 1.3 (v)]{wu}}.

Let $G$ be the subgroup of $PGL(3,q)$ fixing $\mathcal{O}$ setwise.
Then $G$ is the three-dimensional projective orthogonal group over
$\Ff_q$, and it is well known \cite[p.~158]{hir} that $G\cong
PGL(2,q)$. Also, $G$ has an index 2 subgroup $H$, which is
isomorphic to $PSL(2,q)$. It is known \cite{dye} that $H$ acts
transitively on $E$ (respectively, $I$), as well as on $Se$
(respectively, $T$ and $Pa$).

Let $F$ be an algebraic closure of $\Ff_2$. The action of $H$ on $E$ 
makes the vector space $F^{E}$ into an $FH$-permutation module. 
Define
\begin{equation}
\phi: F^E\rightarrow F^E
\end{equation}
by letting $\phi(\p)=\sum_{\q\in \p^{\perp}\cap E}\q$ for each
$\p\in E$, and then extending $\phi$ linearly to $F^E$, where
$\perp$ is the polarity induced by the quadratic form $Q$. Then up
to permutations of the rows and columns, $\A_{33}$ (with its entries
viewed as elements in $F$) is the matrix of $\phi$ with respect to
the basis $E$ of $F^E$. One can easily show that $\phi$ is an
$FH$-homomorphism. Therefore, the $F$-null space of $\A_{33}$ is
equal to $Ker(\phi)$, which is an $FH$-submodule of $F^E$. This
point of view allows us to bring powerful tools from modular
representation theory to bear on problems such as
Conjecture~\ref{con1}. In fact we will not only find the
$F$-dimension of $Ker(\phi)$, but also the $FH$-module structure of
$Ker(\phi)$.

We give a brief overview of the paper. In Section 2, we first prove
several important geometric results related to a conic in $PG(2,q)$,
which allow us to show that if we arrange the rows and columns of
$\A_{33}$ in a particular way, then $\A_{33}^5\equiv \A_{33}\pmod 2$
(Theorem~\ref{matrix}). This latter result has two important
consequences:

(i) the $F$-null space of $\A_{33}$ is equal to the span of the rows
of $\A_{33}^4+I$ (mod 2), where $I$ is the identity matrix,

(ii) as $FH$-modules,
\begin{equation}\label{dirsum}
F^E\cong Ker(\phi)\oplus Im(\phi).
\end{equation}
To prove Conjecture~\ref{con1}, it suffices to compute the dimension
of $Ker(\phi)$. 
In order to do this, we apply Brauer's theory of blocks. The decomposition of
the characters of $H$ into blocks was given by Burkhardt {\cite{burkhardt}}
and Landrock {\cite{landrock}.
We begin by computing the character of the complex permutation module
$\Cc^E$, and its decomposition into blocks. 
This information can be read off from the complex
character table and information (Lemma~\ref{stabilizer_intersection}) about the intersections
of conjugacy classes of $H$ with the subgroup $K$ which stabilizes
an element of $E$. From this we see that $\Cc^E$ is a direct sum
of modules consisting of one simple module from each block of defect zero,
and some summands from blocks of positive defect. 
Then we consider the decomposition of the $Ker(\phi)$ and $Im(\phi)$ into blocks.
According to Brauer's theory, every $FH$-module $M$ is the direct sum
\begin{equation}
M=\displaystyle\bigoplus_{B}e_B M
\end{equation}
where $e_B$ is a primitive idempotent in the center of $FH$.
The block idempotents $e_B$  can be computed as elements of $FH$
from the complex character table of $H$ and the known partition
of the complex characters into blocks. In order to compute $e_BKer(\phi)$
and  $e_BIm(\phi)$
we need detailed information concerning the action of group elements
in various conjugacy classes on various geometric objects
and on the intersections of certain special subsets of
$H$ (see Definition~\ref{definition})  with various conjugacy classes of
$H$. These computations are made in Sections 4 and  5.
This information tells us which block idempotents annihilate
$Ker(\phi)$ and $Im(\phi)$ (Lemma~\ref{y4}). From this, we see that
$Ker(\phi)$ is equal to the direct sum of all of the components
$e_BF^E$ corresponding to blocks of defect zero, or this sum plus 
an additional trivial summand, depending on $q$.
The dimension of $Ker(\phi)$ can then be deduced from the
block decomposition of $\Cc^E$, since the $B$-component of $F^E$ is, 
in a sense which will be made precise, the mod $p$ reduction of the 
$B$-component of $\Cc^E$.
\section{Geometric Results}

In the rest of this paper, to simplify notation, we will use $\B$
(instead of $\A_{33}$) to denote the following $(0,1)$-incidence
matrix between $Se$ and $E$: The columns of $\B$ are labeled by 
the external points $\p_1, \p_2, \ldots ,\p_{q(q+1)/2}$, the rows of
$\B$ are labeled by the secant lines $\p_1^{\perp}, \p_2^{\perp},
\ldots ,\p_{q(q+1)/2}^{\perp}$ and the $(i,j)$-entry of $\B$ is 1 if
and only if $\p_j\in\p_i^{\perp}$. Note that the matrix $\B$ is
symmetric. Our goal in this section is to prove the following
theorem.

\begin{Theorem}\label{matrix}
If we view $\B$ as a matrix over $\Zz$, then
$$\B^5\equiv \B \pmod 2,$$
where the congruence means entrywise congruence. Moreover, if
$q\equiv \pm3 \pmod 8$, then
$$\B^3\equiv \B \pmod 2.$$
\end{Theorem}
\begin{Remark}
As we will see from the proof of Theorem~\ref{matrix}, we do not have
$\B^3\equiv \B \pmod 2$ when $q\equiv \pm1 \pmod 8$.
\end{Remark}

\subsection{Some known geometric results related to a conic in $PG(2,q)$}

In this paper, the classical projective plane $PG(2,q)$ is represented
via homogeneous coordinates. Namely,  a point $\p$ of $PG(2,q)$ can 
be written as $(a_0,a_1,a_2)$, where
$(a_0,a_1,a_2)$ is a nonzero vector of $V$, and a line $\ell$ as
$[b_0,b_1,b_2]$, where $b_0,b_1,b_2$ are not all
zeros. The point $\p=(a_0,a_1,a_2)$ lies
on the line $\ell=[b_0,b_1,b_2]$ if and only if
$$a_0b_0+a_1b_1+a_2b_2=0.$$

Recall that a {\it collineation} of $PG(2,q)$ is an automorphism of
$PG(2,q)$, which is a bijection from the set of all points and all
lines of $PG(2,q)$ to itself that maps a point to a point and a line
to a line, and preserves incidence. It is well known that each
element of $GL(3,q)$, the group of all $3\times 3$ non-singular
matrices over $\Ff_q$, induces a collineation of $PG(2,q)$. The
proof of the following lemma is straightforward.

\begin{Lemma}\label{action1} Let
$\p=(a_0,a_1,a_2)$ $($respectively,
$\ell=[b_0,b_1,b_2])$ be a point $($respectively, a
line$)$ of $PG(2,q)$. Suppose that $\theta$ is a collineation of
$PG(2,q)$ that is induced by $\mathbf{D}\in GL(3,q)$. If we use 
$\p^\theta$ and $\ell^\theta$ to denote the images of $\p$ and 
$\ell$ under $\theta$, respectively, then
$$\p^\theta = (a_0,a_1,a_2)^\theta = (a_0,a_1,a_2)\mathbf{D}$$
and
$$\ell^\theta = [b_0,b_1,b_2]^\theta =[c_0,c_1,c_2],$$
where $c_0,c_1,c_2$ correspond to the first, the second, and 
the third coordinate of the vector $\mathbf{D}^{-1}(b_0,b_1,b_2)^\top$, 
respectively.
\end{Lemma}

A {\it correlation} of $PG(2,q)$ is a bijection from the set of
points to the set of lines as well as the set of lines to the set
of points that reverses inclusion. A {\it polarity}
of $PG(2,q)$ is a correlation of order $2$. The image of a point
$\p$ under a correlation $\sigma$ is denoted by $\p^\sigma$, and
that of a line $\ell$ is denoted by $\ell^\sigma$. It can be shown
\cite[p.~181]{hir} that the non-degenerate quadratic form
$Q(X_0,X_1,X_2)$ = $X_1^2-X_0X_2$ induces a polarity $\sigma$ 
(or $\perp$) of $PG(2,q)$, which can be represented by the matrix
\begin{equation}\label{matrix_M}
\M=\left(\begin{array}{ccc}
0 & 0 & -\frac{1}{2} \\
0 & 1 & 0 \\
-\frac{1}{2} & 0 & 0 \\
\end{array}\right).
\end{equation}

\begin{Lemma}\label{lemB}$(${\rm{\cite[p.~47]{hp}}}$)$
Let $\p=(a_0,a_1,a_2)$ $($respectively,
$\ell=[b_0,b_1,b_2])$ be a point $($respectively, a
line$)$ of $PG(2,q)$. If $\sigma$ is the polarity represented by the
above non-singular symmetric matrix $\M$, then
$$\p^{\sigma} = (a_0,a_1,a_2)^{\sigma} =[c_0,c_1,c_2]$$
and
$$\ell^{\sigma} = [b_0,b_1,b_2]^{\sigma} =(b_0,b_1,b_2)\M^{-1},$$
where $c_0,c_1,c_2$ correspond to the first, the second, the third 
coordinate of the column vector $\M(a_0,a_1,a_2)^\top$, respectively.
\end{Lemma}
For example, if $\p=(x,y,z)$ is a point of $PG(2,q)$,
then its image under $\sigma$ is $\p^{\sigma}=[z,-2y,x]$.

For convenience, we will denote the set of all non-zero squares of
$\Ff_q$ by $\Sq_q$, and the set of non-squares by $\Nsq_q$. Also,
$\Ff_q^*$ is the set of non-zero elements of $\Ff_q$.

\begin{Lemma}\label{bijection}$(${\rm{\cite[p.~181--182]{hir}}}$)$
Assume that $q$ is odd.
\begin{itemize}
\item[(i)] The polarity $\sigma$ above defines three bijections; that is,
$\sigma:\;I\rightarrow\;Pa$, $\sigma:\;E\rightarrow\;Se$, and
$\sigma:\;\mathcal{O}\rightarrow\;T$ are all bijections.
\item[(ii)] A line $[b_0,b_1,b_2]$ of $PG(2,q)$ is a passant,
a tangent, or a secant to $\mathcal{O}$ if and only if 
$b_1^2-4b_0b_2 \in\Nsq_q$, $b_1^2-4b_0b_2 = 0$, or 
$b_1^2-4b_0b_2\in \Sq_q$, respectively.
\item[(iii)] A point $(a_0,a_1,a_2)$ of $PG(2,q)$ is internal,
absolute, or external if and only if $a_1^2-a_0a_2 \in \Nsq_q$, 
$a_1^2-a_0a_2 =0$, or $a_1^2-a_0a_2 \in \Sq_q$, respectively.
\end{itemize}
\end{Lemma}

The results in the following lemma can be obtained by simple
counting; see {\rm \cite{hir}} for more details and related results.
\begin{Lemma}{\rm (\cite[p.~170]{hir})}
Using the above notation, we have
\begin{equation}\label{number}
|T|=|\mathcal{O}|=q+1,\;|Pa|=|I|=\frac{q(q-1)}{2},\;\text{and}\;|Se|=|E|=\frac{q(q+1)}{2}.
\end{equation}
Also, we have the following tables:
\begin{table}[htp]
\begin{center}
\caption{Number of points on lines of various types}
\bigskip
\begin{tabular}{cccc}
\hline
{Name} & {Absolute points} & {External points} & {Internal points} \\
\hline
{Tangent lines} & $1$ & $q$ & $0$ \\
{Secant lines} & $2$ & $\frac{q-1}{2}$ & $\frac{q-1}{2}$ \\
{Passant lines} & $0$ & $\frac{q+1}{2}$ & $\frac{q+1}{2}$\\
\hline
\end{tabular}
\label{tab1}
\end{center}
\end{table}

\begin{table}[htp]
\begin{center}
\caption{Number of lines through points of various types}
\bigskip
\begin{tabular}{cccc}
\hline
{Name} & {Tangent lines} & {Secant lines} & {Skew lines} \\
\hline
{Absolute points} & $1$ & $q$ & $0$ \\
{External points} & $2$ & $\frac{q-1}{2}$ & $\frac{q-1}{2}$ \\
{Internal points} & $0$ & $\frac{q+1}{2}$ & $\frac{q+1}{2}$\\
\hline
\end{tabular}
\label{tab2}
\end{center}
\end{table}
\end{Lemma}

\subsection{More geometric results}

Let $G$ be the automorphism group of $\mathcal{O}$ in $PGL(3,q)$
(i.e. the subgroup of $PGL(3,q)$ fixing ${\mathcal O}$ setwise). 
We define
\begin{equation}\label{groupG}
\mathcal{G}:=
\left.\left\{\left(\begin{array}{ccc}
a^2 & ab & b^2\\
2ac & ad+bc & 2bd \\
c^2& cd & d^2\\
 \end{array}\right)\right| a, b, c, d \in \Ff_q, ad-bc\in\Ff_q^*\right\}.
\end{equation}

For the convenience of the reader, we include the discussion of the 
structure of $G$ which can be found in \cite[p.158]{hir}. The 
conic $\mathcal{O}=\{(1,t,t^2)\mid t\in \Ff_q\}\cup\{(0,0,1)\}$ can be 
identified with the projective line $PG(1,q)$ via $(1,t)\leftrightarrow(1,t,t^2)$
and $(0,1)\leftrightarrow (0,0,1)$, where the point $(0,1)$ of $PG(1,q)$, 
usually denoted by $\infty$, is the ``point at infinity" of $PG(1,q)$. 
A generic element $\varphi$ of $PGL(2,q)$ is represented by 
$(1,t)\mapsto (1,t) {\bf M}(\varphi)$, where 
${\bf M}(\varphi)=\left(\begin{smallmatrix}a & b \\ c & d\end{smallmatrix}\right)$
with $ad-bc\not=0$. This mapping may be written as $t\mapsto\frac{dt+b}{ct+a}$
with the usual conventions on $\infty$. It has the following effect on $\mathcal{O}$:
$$(1,t,t^2)\mapsto ((ct+a)^2, (dt+b)(ct+a), (dy+b)^2)=(1,t,t^2)^{\varphi^{'}},$$
where $\phi^{'}$ is represented by
\begin{displaymath}
{\bf M}(\varphi^{'})=
\left(\begin{array}{ccc}
a^2 & ab & b^2\\
2ac & ad+bc & 2bd \\
c^2& cd & d^2\\
 \end{array}\right).
\end{displaymath}
So $\varphi$ induces an element $\varphi^{'}$ of $G$. An easy counting 
argument can be applied to show that there are $q^5-q^2$ conics in 
$PG(2,q)$. Hence, $|G|=\frac{|PGL(3,q)|}{q^5-q^2}=|PGL(2,q)|$. So 
$\tau: PGL(2,q)\rightarrow G$ given by $\tau(\varphi)=\varphi^{'}$ is a bijection. 
It is straightforward to check that $\tau$ is a 
homomorphism. We have shown the following lemma.
\begin{Lemma}
The above $\mathcal{G}$ induces $G$ and $G\cong PGL(2,q)$.
\end{Lemma}

From the above discussions, we have an isomorphism $\tau$ of $PGL(2,q)$
to the group $G$ induced by $\mathcal{G}$. From now on, we will identify $G$ and 
$\mathcal{G}$. Recall that $PSL(2,q)=SL(2,q)/\{\pm\left(\begin{smallmatrix}1& 0\\ 0 & 1
\end{smallmatrix}\right)\}$.

In the rest of the paper, we always use $\xi$ to denote a fixed primitive 
element of $\Ff_q$. For $a,b,c\in\Ff_q$, we define
\begin{equation*}
{\bf d}(a,b,c):=
\begin{pmatrix}
a & 0 & 0 \\
0 & b & 0 \\
0 & 0 & c\\
\end{pmatrix},\;
{\bf ad}(a,b,c):=
\begin{pmatrix}
0& 0 & a\\
0 & b & 0\\
c & 0 & 0\\
\end{pmatrix}.
\end{equation*}

Next we define
\begin{equation}\label{grouph}
H:=
\left.\left\{\left(\begin{array}{ccc}
a^2 & ab & b^2\\
2ac & ad+bc & 2bd \\
c^2& cd & d^2\\
 \end{array}\right)\right| a, b, c, d \in \Ff_q, ad-bc=1\right\}\subset \mathcal{G}.
\end{equation}
It is obvious that the image of $PSL(2,q)$ under $\tau$ is $H$, and thus 
$H\cong PSL(2,q)$. The isomorphism between
$H$ and $PSL(2,q)$ was proved by Dickson (\cite[Theorem 178]{dickson})
in different coordinates.

Since $$PGL(2,q)=PSL(2,q)\cup 
\begin{pmatrix}1& 0\\ 0 & \xi^{-1}\end{pmatrix}\cdot PSL(2,q),$$
by the isomorphism $\tau$ we have
\begin{equation}\label{groupg}
G=H\cup {\bf d}(1,\xi^{-1},\xi^{-2})\cdot H.
\end{equation}
Moreover, the following holds.
\begin{Lemma}\label{transitive}\cite{dye}
The group $G$ acts transitively on both $I$
$($respectively, $Pa)$ and $E$ $($respectively, $Se)$.
\end{Lemma}

\begin{Lemma}\label{meet}
Let $\p$ be a point not on ${\mathcal O}$, $\ell$ be a non-tangent
line, and $\p\in\ell$. Using the above notation, we have the
following.

\begin{enumerate}
\renewcommand{\labelenumi}{(\roman{enumi})}

\item If $\p\in I$ and $\ell\in$Pa, then $\p^\perp\cap \ell \in E$
if $q\equiv 1 \pmod 4$, and $\p^\perp \cap \ell \in I$ if $q\equiv 3
\pmod 4$.

\item If $\p\in I$ and $\ell \in$Se, then $\p^\perp\cap \ell \in I$
if $q\equiv 1 \pmod 4$, and $\p^\perp\cap \ell \in E$ if $q\equiv 3
\pmod 4$.

\item If $\p\in E$ and $\ell \in$Pa, then $\p^\perp\cap \ell \in I$
if $q\equiv 1 \pmod 4$, and $\p^\perp\cap \ell \in E$ if $q\equiv 3
\pmod 4$.

\item If $\p\in E$ and $\ell \in$Se, then $\p^\perp\cap \ell\in E$
if $q\equiv 1 \pmod 4$, and $\p^\perp\cap\ell\in I$ if $q \equiv 3
\pmod 4$.

\end{enumerate}
\end{Lemma}
{\Proof} The proofs of all four parts are similar. So we only give
the proof of part (i). By Lemma~\ref{transitive}, $G$ acts
transitively on $I$. Thus we may assume that
$\p=(1,0,c)$, where $-c\in\Nsq_q$. By assumption
$\p\in \ell$ and $\ell\in Pa$, we have $\ell =
[1,b,-c^{-1}]$ for some $b\in\Ff_q$ such that
$\frac{4+b^2 c}{c}\in \Nsq_q$.

Now $\p^\perp=[-\frac{1}{2}c, 0, -\frac{1}{2}]$ by Lemma~\ref{lemB}. 
So $\p^\perp\cap \ell =(b, -2, -bc)$. If $q\equiv 1\pmod 4$, then
$-1\in\Sq_q$. From $-c\in\Nsq_q$ we now see that $c\in \Nsq_q$.
Thus, in this case, the condition $\frac{4+b^2 c}{c} \in \Nsq_q$
implies that $4+b^2 c\in \Sq_q$. By Lemma~\ref{bijection} (iii), we
have $\p^\perp\cap \ell \in E$. If $q\equiv 3 \pmod 4$, then $-1\in
\Nsq_q$. Similar arguments show that $\p^\perp\cap\ell \in I$. \QED

We define $\Sq_q-1:=\{s-1\mid s\in \Sq_q\}$
and $\Nsq_q-1:=\{s-1\mid s\in \Nsq_q\}$.  
In the rest of the paper, we will frequently use the following lemma.
\begin{Lemma}\label{cs}\cite{store}
Using the above notation,
\begin{itemize}
\item[(i)] if $q\equiv 1\pmod 4$, then $|(\Sq_q-1)\cap\Sq_q|=\frac{q-5}{4}$
and $|(\Sq_q-1)\cap\Nsq_q|$ $=$ $|(\Nsq_q-1)\cap\Sq_q|$ $=$ $|(\Nsq_q-1)\cap\Nsq_q|
=\frac{q-1}{4}$;
\item[(ii)] if $q\equiv 3\pmod 4$, then $|(\Nsq_q-1)\cap\Sq_q|=\frac{q+1}{4}$ and
$|(\Sq_q-1)\cap\Sq_q|$ $=$ $|(\Sq_q-1)\cap\Nsq_q|$ $=$ $|(\Nsq_q-1)\cap\Nsq_q|
=\frac{q-3}{4}$.
\end{itemize}
\end{Lemma}

The following lemma will be used throughout the rest of this section
and Section 3. Its proof is straightforward. We omit the details.
\begin{Lemma}\label{stabilizer}
Let $W$ be a subgroup of $G$. Suppose that $g\in G$ and $\p$ is a
point of $PG(2,q)$. Then
$$Stab_{W^g}(\p^g) = Stab_W(\p)^g,$$
where $Stab_W(\p)^g=\{g^{-1}hg\mid h \in Stab_W(\p)\}$ and $W^g=\{g^{-1}
hg\mid h \in W\}$.
\end{Lemma}

Let $\p$ be a point of $PG(2,q)$. Then $\p^{\perp}$ is a line of
$PG(2,q)$. We will use $I_{\p^{\perp}}$ (respectively,
$E_{\p^{\perp}}$ and $\mathcal{O}_{\p^\perp}$) to denote the set of internal 
(respectively, external and conic) points on $\p^{\perp}$. Also we will use $Pa_{\p}$
(respectively, $Se_{\p}$ and $T_{\p}$) to denote the set of passant
(respectively, secant and tangent) lines through $\p$. Then the following
lemma is apparent.

\begin{Lemma}\label{one-to-one} Let $\p$ be a point of $PG(2,q)$. Then
the polarity $\perp$ defines a bijection between
 $I_{\p^\perp}$ and $Pa_{\p}$,  and it also defines a bijection between 
 $E_{\p^\perp}$ and $Se_\p$.
\end{Lemma}
From simple matrix computations, we have the following lemma.
\begin{Lemma}\label{com}
Let $\p$ be a point of $PG(2,q)$, $\perp$ the polarity of $PG(2,q)$
defined above, and $g\in G$. Then $(\p^\perp)^g=(\p^g)^\perp$.
\end{Lemma}

\begin{Lemma}\label{orderK}
Let $\p\in E$, $K=Stab_G(\p)$ $($i.e. the stabilizer of $\p$ in G$)$, and
$\p_1\in \p^\perp$.
\begin{itemize}
\item[(i)] If $\p_1$ is an internal or external point on $\p^\perp$,
then $|Stab_K(\p_1)|=4$.
\item[(ii)] If $\p_1$ is a point of $\mathcal{O}$, then $|Stab_K(\p_1)| = q-1$.
\end{itemize}
\end{Lemma}
{\Proof} Let $\q=(0,1,0)\in E$, and $K'=Stab_G(\q)$.
By Lemma~\ref{transitive}, we can find $g\in G$ such that $\q=\p^g$.
Now let
$\p_1\in \p^\perp$. Then $\q_1:=\p_1^g\in \q^\perp$. By
Lemma~\ref{stabilizer}, we have
$$Stab_{K^{'}}(\q_1) = Stab_{K^g}(\p_1^g) = Stab_{K}(\p_1)^g.$$
It follows that $|Stab_K(\p_1)| = |Stab_{K^{'}}(\q_1)|$.
Therefore, in order to prove the lemma, it is enough to calculate
$Stab_{K^{'}}(\q_1)$, where $\q_1$ is an arbitrary point on $\q^\perp$.

First we calculate $K^{'}=Stab_G(\q)$. By (\ref{groupg}), we have
 $G=H\cup {\bf d}(1,\xi^{-1},\xi^{-2})H$, where $H$
is defined in (\ref{grouph}).
Let
\begin{displaymath}
g=\left(\begin{array}{lll}
a^2 & ab & b^2 \\
2ac & ad+bc & 2bd \\
c^2 & cd & d^2 \\
\end{array}\right) \in H
\end{displaymath}
Then by Lemma~\ref{action1}, we see that $g$ fixes $\q$ if and only if
\begin{equation*}\label{eq3.1}
( 0, 1, 0)g= u(0,1,0)
\end{equation*}
for some $u\in\Ff_q^*$.
Since $ad-bc=1$, we have either $b=c=0$, $ad=1$
or $a=d=0$, $bc=-1$. If $g\in {\bf d}(1,\xi^{-1},\xi^{-2})H$, then similar 
calculations yield the same two cases. Therefore,
\begin{equation}\label{stab}
\begin{array}{lll}
K^{'} &= &\left\{ {\bf ad}(c^{-2}, -1, c^2), {\bf ad}(c^{-2},
-\xi^{-1}, c^2\xi^{-2})\mid c \in \Ff_q^*\right\}\\
{} & {} &\cup \left\{{\bf d}(d^2, 1, d^{-2}), {\bf d}(d^2, \xi^{-1}, d^{-2}\xi^{-2})\mid d \in\Ff_q^*\right\}.
\end{array}
\end{equation}

It is clear that $\q^{\perp}=[0,1,0]$. Thus
$I_{\q^\perp}=\{(1,0,-n)\mid n\in \Nsq_q\}$ and
$E_{\q^\perp}=\{(1,0,-m)\mid m \in \Sq_q\}$. 

Let
$\q_1=(1,0,-n)\in I_{\q^\perp}$ (so $n\in\Nsq_q$).
Then
\begin{equation*}
\begin{array}{lll}
Stab_{K^{'}}(\q_1) & = & \left\{{\bf d}(1,1,1), {\bf d}(-1,1,-1)\right\}\\
{} & {}&\cup \left\{{\bf ad}(n\xi^{-1},-\xi^{-1},(n\xi)^{-1}), {\bf ad}(-n\xi^{-1},
-\xi^{-1},-(n\xi)^{-1}) \right\}
\end{array}
\end{equation*}
if $q\equiv 1\pmod 4$, and 
\begin{equation*}
\begin{array}{lll}
Stab_{K^{'}}(\q_1) & = & \left\{ {\bf d}(1,1,1), {\bf d}(-\xi^{-1},
\xi^{-1},-\xi^{-1})\right\}\\
{} & {} & \cup \left\{{\bf ad}(-n,-1,-n^{-1}), {\bf
ad}(n\xi^{-1},-\xi^{-1},(n\xi)^{-1}) \right\}.
\end{array}
\end{equation*}
Hence, $|Stab_{K^{'}}(\q_1)| = 4$ if $\q_1\in I_{\q^\perp}$.

Let $\q_1 =(1,0,-m)\in E_{\q^\perp}$ (so $m\in
\Sq_q$). Then
\begin{equation*}
\begin{array}{lll}
Stab_{K^{'}}(\q_1)& = & \left\{ {\bf d}(1,1,1), {\bf d}(-1,1,-1), {\bf
ad}(-m,-1,-m^{-1}), {\bf ad}(m,-1,m^{-1}) \right\}
\end{array}
\end{equation*}
if $q\equiv 1\pmod 4$, and
\begin{equation*}
\begin{array}{lll}
Stab_{K^{'}}(\q_1) 
&=& \left\{ {\bf d}(1,1,1), {\bf ad}(-m,-1,-m^{-1})\right\}\\
{} & {} &\cup \left\{{\bf d}(-\xi^{-1},\xi^{-1},-\xi^{-1}), {\bf ad}(-m\xi^{-1},-\xi^{-1},-(m\xi)^{-1}) \right\}
\end{array}
\end{equation*}
if $q\equiv 3\pmod 4$.
Hence, $|Stab_{K^{'}}(\q_1)| =4$ if $\q_1\in E_{\q^\perp}$.

The proof of $(ii)$ is similar. We omit the details. \QED

\begin{Corollary}\label{stab_1}
Let $\p\in E$. Then $|Stab_G(\p)|= 2(q-1)$.
\end{Corollary}
{\Proof} The corollary follows from the fact that $|K^{'}|=2(q-1)$, where
$K^{'}$ is the group given in (\ref{stab}). 
\QED

\begin{Corollary}\label{commutative}
Let $\p$ be a point of $PG(2,q)$ and let $\perp$ be the
polarity of $PG(2,q)$ defined above. Then for $g\in Stab_G(\p)$ we have
$\p^\perp = (\p^\perp)^g$. Consequently, $\p^\perp$ is fixed setwise by
$Stab_G(\p)$. Moreover, $Stab_G(\p^\perp)=Stab_G(\p)$.
\end{Corollary}

{\Proof } The first part of the corollary follows from Lemma~\ref{com}. Therefore,
$Stab_G(\p)\subseteq Stab_G(\p^\perp)$. Assume that $\p\in E$. Since $G$
acts transitively on $Se$ and $E$ respectively, we have
$$|Stab_G(\p^\perp)| = \frac{|G|}{|Se|} = 2(q-1) = |Stab_G(\p)|,$$
where the last equality follows from Corollary~\ref{stab_1}.
Thus, $Stab_G(\p)= Stab_G(\p^\perp)$. The case where $\p\in I$
or $\p\in\mathcal{O}$ can be proved in the same way.
\QED

\begin{Proposition}\label{K_transitive}
Let $\p\in E$ and $K=Stab_G(\p)$. Then $K$ is transitive on $I_{\p^\perp}$,
$E_{\p^\perp}$ and $\mathcal{O}_{\p^\perp}$. Also $K$ is transitive on 
$Pa_{\p}$, $Se_{\p}$ and $T_{\p}$.
\end{Proposition}

{\Proof} Let $\p_1$ be an internal or external point on $\p^{\perp}$.
By Lemma~\ref{orderK}, the length of the orbit of $\p_1\in \p^\perp$
under the action of $K$ is
$$\frac{|K|}{|Stab_K(\p_1)|} =\frac{2(q-1)}{4}=\frac{q-1}{2},$$
which is equal to $|I_{\p^\perp}|=|E_{\p^\perp}|$. So $K$ is
transitive on both $I_{\p^\perp}$ and $E_{\p^\perp}$. From
Corollary~\ref {commutative}, it follows that $K$ is also transitive
on both $Pa_{\p}$ and $Se_{\p}$. The case where $\p_1\in\mathcal{O}_{\p^\perp}$
can be proved in the same way.\QED

Let $\p\in E$, $\ell\in Se$, and $\p\not\in\ell$. We use $Se_E(\p,\ell)$ to denote
the set of secant lines through $\p$ meeting $\ell$ in an external point. That is,
$$Se_E(\p,\ell)=\{\ell_1\in Se_\p\mid \ell_1\cap\ell\in E\}.$$

\begin{Lemma}\label{intersection}
Let $\ell$ be a secant line and $T_1$, $T_2$ the two tangent lines through
$\p^{'}:=\ell^\perp$. Suppose that $\p\in E$, $\p\notin \ell$, and $\p \ne \p^{'}$.
\begin{itemize}
\item[(i)] If $\p$ is on either $T_1$ or $T_2$, then $|Se_E(\p,\ell)|$ is odd  or 
even according as $q\equiv \pm 1\pmod 8$ or $q\equiv \pm 3\pmod 8$.
\item[(ii)] If $\p$ is on a passant or a secant through $\p^{'}$, then
$|Se_E(\p,\ell)|$ is even.
\end{itemize}
\end{Lemma}

{\Proof} Since $G$ acts transitively on $Se$ and preserves incidence, 
we may take $\ell=[0,1,0]$. Then $\p^{'}=\ell^\perp=(0,1,0)$,
$T_1=[0,0,1]$, and $T_2=[1,0,0]$.
Also, $$Se_{\p^{'}}=\{[1,0,y]\mid y\in \Ff_q^*,\;-4y\in\Sq_q\}$$
and
$$Pa_{\p^{'}}=\{[1,0,x]\mid x\in\Ff_q^*,\;-4x\in\Nsq_q\}.$$

Since $K=Stab_G(\p^{'})$ acts transitively on $Pa_{\p^{'}}$,
$Se_{\p^{'}}$, and $\{T_1,T_2\}$ by Proposition~\ref{K_transitive},
 we see that, in order to prove the lemma, it
is enough to consider the external points excluding $\p^{'}$ on a
{\it special} secant, passant, and tangent line through $\p^{'}$. To
this end, we take $\ell_1=[1,0,y]$,
$y\in\Ff_q^*$, $-4y\in\Sq_q$, $\ell_2=[1,0,x]$,
$x\in\Ff_q^*$, $-4x\in\Nsq_q$, and $\ell_3=T_2$ to be the special
secant, passant, and tangent line through $\p'$, respectively. The
external points excluding $\p^{'}$ on $\ell_1$, $\ell_2$, and
$\ell_3$ but not on $\ell$ are given, respectively, by
\begin{equation*}\label{spoints}
\begin{array}{ll}
E_{\ell_1}=\{(1,m,-y^{-1})\mid m\in \Ff_q^*,\;-4y\in\Sq_q,\;m^2+y^{-1}\in\Sq_q\},\\
E_{\ell_2}=\{(1,n,-x^{-1})\mid n \in \Ff_q^*,\;-4x\in\Nsq_q,\;m^2+x^{-1}\in\Sq_q\},\\
E_{\ell_3}=\{(0,1,s)\mid s\in \Ff_q^*\}.\\
\end{array}
\end{equation*}
To prove the lemma, we may assume that $\p$ is in $E_{\ell_1}$, $E_{\ell_2}$, or
$E_{\ell_3}$.

If $\p=(1,m,-y^{-1}) \in E_{\ell_1}$, then by (\ref{stab}), we obtain that
\begin{equation*}
Stab_K(\p)=\left\{ {\bf d}(1,1,1), {\bf ad}(y^{-1}, -1, y) \right\}
\end{equation*}
if $q\equiv 1\pmod4$, and
\begin{equation*}
Stab_K(\p)=\left\{ {\bf d}(1, 1, 1), {\bf ad}((y\xi)^{-1}, -\xi^{-1},
y\xi^{-1}) \right\}
\end{equation*}
if $q\equiv 3\pmod 4$. In particular, $|Stab_K(\p)|=2.$

The lines through $\p$ are
$$\{[1,n_1,y(1+mn_1)]\mid n_1\in \Ff_q\}\cup\{[0,1,my]\}.$$

From
\begin{displaymath}
{\bf ad}(y^{-1}, -1, y)^{-1}(1,n_1,y(1+mn_1))^\top = (1+mn_1, -n_1,
y)^\top
\end{displaymath}
if $q\equiv 1\pmod 4$
and
\begin{displaymath}
{\bf ad}((y\xi)^{-1}, -\xi^{-1}, y\xi^{-1})^{-1}
(1,n_1,y(1+mn_1))^\top = (\xi(1+mn_1), -n_1\xi, y\xi)^\top
\end{displaymath}
if $q\equiv 3\pmod 4$,
it follows that a line of the form $[1,n_1,y(1+mn_1)]$ 
is fixed by $Stab_K(\p)$ if and only if
\begin{displaymath}
\left\{\begin{array}{lll}
\frac{-n_1}{1+mn_1} &= & n_1,\\
\frac{y}{1+mn_1} & = &y(1+mn_1).\\
\end{array}\right.
\end{displaymath}
From the two equations, we obtain that $n_1=-2m^{-1}$. Therefore 
$\ell^{'}:=[1,-2m^{-1},-y]$
is the unique line of the form $[1,n_1,y(1+mn_1)]$ 
through $\p$ that is fixed by $Stab_K(\p)$. Easy calculations now show 
that $[0,1,my]$ cannot be fixed by $Stab_K(\p)$. 
Also note that $\ell_1$ is fixed by $Stab_K(\p)$. Thus, under the action 
of $Stab_K(\p)$, the lines through $\p$ are split into $\frac{q+3}{2}$ orbits, 
two of which have length $1$ (namely, $\{\ell_1\}$ and $\{\ell'\}$),
and $\frac{q-1}{2}$ of which have length $2$. Lines in the same orbit of 
length $2$ must be of the same type; that is, they must be both secants, 
or both passants, or both tangents.

When $q\equiv 1 \pmod 4$, we have $\ell^{'}\in Se_{\p}$, $\ell_1\in Se_{\p^{'}}$,
$\ell_1\cap \ell=(1,0,-y^{-1})\in E$, and 
$\ell^{'}\cap\ell=(1,0,y^{-1})\in E$ since $-1\in\Sq_q$, $y\in \Sq_q$ 
and $m^2+y^{-1}\in \Sq_q$. In this case, $|Se_E(\p,\ell)|$ is even.

When $q\equiv 3\pmod 4$, we have $\ell_1\in Se_{\p^{'}}$, $\ell^{'}\in Pa_{\p}$,
$\ell_1\cap\ell\in I$, and $\ell^{'}\cap\ell\in E$ since  $-1\in\Nsq_q$, $y\in \Nsq_q$ 
and $m^2+y^{-1}\in \Sq_q$. Therefore, $|Se_E(\p,\ell)|$ is even as well.

The proof of the lemma in the case where $\p\in\ell_1$ is now finished. Similar 
arguments can be applied to show that $|Se_E(\p,\ell)|$ is even if $\p\in\ell_2$. 
We omit the details.

The rest of the proof is concerned with the parity of $|Se_E(\p,\ell)|$ when 
$\p\in \ell_3$. Let $\p=(0,1,s)\in \ell_3$ with $s\not=0$. The set of lines 
connecting $\p$ with external points on $\ell$ is
$$L_s=\left\{[1,su^{-1},-u^{-1}]\mid -u\in \Sq_q\right\}.$$
The number of secant lines in $L_s$ is determined by the number of
$u$ satisfying both of the following two conditions
\begin{displaymath}
\begin{array}{ll}
(a) & \frac{s^2+4u}{u^2} \in \Sq_q,\\
(b) & -u\in \Sq_q.\\
\end{array}
\end{displaymath}

If $q\equiv 1 \pmod 4$, the number of $u$ satisfying $(a)$ and $(b)$
is equal to $$|\Sq_q\cap (\Sq_q-1)|=\frac{q-5}{4}$$ by Lemma~\ref{cs} (i).

If $q\equiv 3 \pmod 4$, the number of $u$ satisfying $(a)$ and $(b)$
is equal to $$|\Nsq_q\cap (\Sq_q-1)|=\frac{q-3}{4}$$ by Lemma~\ref{cs} (ii).

Therefore, when $q\equiv \pm 1 \pmod 8$, $|Se_E(\p,\ell)|$ is odd;
when $q\equiv \pm 3 \pmod 8$, $|Se_E(\p,\ell)|$ is even. The proof 
is now complete.
\QED

\begin{Definition}
Let $\p\in E$. We define
\begin{displaymath}
N_E(\p)=
\begin{cases}
\{\q\in E\mid \q \in \ell,\;\ell\in Se_{\p}\}\setminus\{\p\}, &\text{if}\;q\equiv 1\pmod 4,\\
\{\q\in E\mid \q \in \ell,\;\ell\in Se_{\p}\}, &\text{if}\;q\equiv 3\pmod 4.\\
\end{cases}
\end{displaymath}
That is, $N_E(\p)$ is the set of external points on the secant lines
through $\p$, from which $\p$ is excluded or not depending on
whether $q\equiv 1\pmod 4$ or $q\equiv 3\pmod 4$. Informally, the
set $N_E(\p)$ can be thought of as the set of {\it external
neighbors} of $\p$.
\end{Definition}

In the following lemma, we investigate the parity of the intersection of 
the external neighbors of two distinct external points, which plays a 
crucial role in the proof of Theorem~\ref{matrix}.

\begin{Lemma}\label{intersection2}
Let $\p_1$ and $\p_2$ be two distinct external points and $\ell_{\p_1,\p_2}$ 
the line through $\p_1$ and $\p_2$.
\begin{itemize}
\item[(i)] If $\ell_{\p_1,\p_2}\in Pa$, then $|N_E(\p_1)\cap N_E(\p_2)|$ is even.
\item[(ii)] If $\ell_{\p_1,\p_2}\in Se$, then $|N_E(\p_1)\cap N_E(\p_2)|$ is odd.
\item[(iii)] If $\ell_{\p_1,\p_2} \in T$, then $|N_E(\p_1)\cap N_E(\p_2)|$ is odd 
or even according as $q\equiv \pm 1\pmod 8$ or $q\equiv \pm 3\pmod 8$.
\end{itemize}
\end{Lemma}

{\Proof} We first give the proof of (ii). Suppose that $\ell_{\p_1,\p_2}\in Se$. 
We consider two subcases.

\noindent{\bf Case I}. $\p_2\notin \p_1^\perp$. 

Let $T_1$, $T_2$ be the two tangent lines through $\p_1$ and
$\mathcal{O}_1=T_1\cap \p_1^\perp = T_1^\perp$,
$\mathcal{O}_2=T_2\cap\p_1^\perp=T_2^\perp$ the two points
of $\mathcal{O}$ on $\p_1^\perp$. Then $\ell_{\p_2,\mathcal{O}_1}$ and
$\ell_{\p_2,\mathcal{O}_2}$ are two secant lines. Also, 
$\ell_{\p_2,\mathcal{O}_1}^\perp\in T_1$
and $\ell_{\p_2,\mathcal{O}_2}^\perp \in T_2$ since
$T_1^\perp=\mathcal{O}_1\in \ell_{\p_2,\mathcal{O}_1}$ and
$T_2^\perp=\mathcal{O}_2\in\ell_{\p_2,\mathcal{O}_2}$. Hence both
$|N_E(\p_1)\cap E_{\ell_{\p_2,\mathcal{O}_1}}|=|Se_E(\p_1,\ell_{\p_2,\mathcal{O}_1})|$ 
and $|N_E(\p_1)\cap E_{\ell_{\p_2,\mathcal{O}_2}}|=|Se_E(\p_1,\ell_{\p_2,\mathcal{O}_2})|$ 
are odd if $q\equiv \pm 1\pmod 8$ and they are even if $q\equiv \pm 3 \pmod 8$
by Lemma~\ref{intersection} (i). Next we consider $\ell_s\in Se_{\p_2}$ such that
$\ell_s\not=\ell_{\p_2,\mathcal{O}_1}$, $\ell_{\p_2,\mathcal{O}_2}$, $\ell_{\p_1,\p_2}$.
Note that the line $\ell_{\p_1,\ell_s^\perp}$ cannot be a tangent line 
since $\ell_s^\perp$ is on neither $T_1$ nor $T_2$ and there are 
only two tangent lines through $\p_1$. Thus $|N_E(\p_1)\cap E_{\ell_s}|$ 
is even by Lemma~\ref{intersection} (ii). It is also clear that
\begin{displaymath}
|N_E(\p_1)\cap E_{\ell_{\p_1,\p_2}}|=
\begin{cases}
\frac{q-3}{2}, & \text{if}\;q\equiv 1\pmod 4, \\
\frac{q-1}{2},  & \text{if}\;q\equiv 3\pmod 4,\\
\end{cases}
\end{displaymath}
which shows that $|N_E(\p_1)\cap E_{\ell_{\p_1,\p_2}}|$ is odd. Set
$L:=Se_{\p_2}\setminus\{\ell_{\p_2,\mathcal{O}_2}, \ell_{\p_2,\mathcal{O}_1},\ell_{\p_1,\p_2}\}$.
Then $|L|=\frac{q-1}{2}-3$ and
\begin{displaymath}
\begin{array}{lllllll}
{} & |N_E(\p_1)\cap N_E(\p_2)|& {}\\
=&
\begin{cases}
\displaystyle\sum_{\ell\in Se_{\p_2}}|(E_{\ell}\setminus\{\p_2\})\cap N_E(\p_1)|, & \text{if}\;
q\equiv 1\pmod 4\\
\displaystyle\sum_{\ell\in Se_{\p_2}}|E_{\ell}\cap N_E(\p_1)| -\frac{q-3}{2}, & \text{if}\;q\equiv 3\pmod 4\\
\end{cases}\\
=&
\begin{cases}
\displaystyle\sum_{\ell\in L}(|E_{\ell}\cap N_E(\p_1)|-1) + (k_1-1)+(k_2-1) +(\frac{q-3}{2}-1), &
 \text{if}\;q\equiv 1 \pmod 4\\
\displaystyle\sum_{\ell\in L}|E_{\ell}\cap N_E(\p_1)|+ k_1+k_2 +\frac{q-1}{2}-\frac{q-3}{2}, &
\text{if}\;q\equiv 3\pmod 4,\\
\end{cases}\\
\end{array}
\end{displaymath}
where $k_1=|E_{\ell_{\p_2,\mathcal{O}_1}}\cap N_E(\p_1)|$ and $k_2=|E_{\ell_{\p_2,\mathcal{O}_2}}\cap N_E(\p_1)|$ are odd if $q\equiv \pm 1\pmod 8$;
otherwise, these numbers are even. Note that $|E_{\ell}\cap N_E(\p_1)|$, $\ell\in L$ are all even by the above discussion. Thus
$$|N_E(\p_1)\cap N_E(\p_2)|\equiv 1\pmod 2.$$

\noindent{\bf Case II.} $\p_2\in \p_1^\perp$. 

This case can happen
only when $q\equiv 1 \pmod 4$ by applying Lemma~\ref{meet} (iv)
to the incidence pair $(\p_1,\ell_{\p_1, \p_2})$. 

Let $\ell_s\in Se_{\p_2}$ such that $\ell_s\not=\ell_{\p_1,\p_2}$ or
$\p_1^\perp$. Then $\ell_s^\perp\notin T_1$ or $T_2$ since,
otherwise, $\ell_s\cap T_1\in \mathcal{O}$ or $\ell_s\cap T_2\in
\mathcal{O}$, which is not the case. The line $\ell_{\p_1,
\ell_s^\perp}$ must be either a secant or a passant. Hence
$|N_E(\p_1)\cap E_{\ell_s}|=|Se_E(\p_1,\ell_s)|$ is even by
Lemma~\ref{intersection} (ii). Note that since $q\equiv 1\pmod 4$,
$|N_E(\p_1)\cap E_{\ell_{\p_1,\p_2}}|= \frac{q-3}{2}$ is odd. If
$\ell_s=\p_1^\perp$, then $|N_E(\p_2)\cap E_{\ell_s}|=\frac{q-1}{2}$
for $q\equiv 1\pmod 4$. Set
$L:=Se_{\p_2}\setminus\{\ell_{\p_1,\p_2}, \p_1^\perp\}$. Then
$|L|=\frac{q-1}{2}-2$ and
\begin{displaymath}
\begin{array}{lll}
{}&|N_E(\p_1)\cap N_E(\p_2)| &{}\\
 = & \displaystyle\sum_{\ell\in L} (|N_E(\p_1)\cap E_{\ell}|-1) + (\frac{q-3}{2}-1) + (\frac{q-1}{2}-1)& {}\\
\equiv & 1 \pmod 2, & {}
\end{array}
\end{displaymath}
where $|N_E(\p_1)\cap E_{\ell}|$, $\ell\in L$, are all even numbers 
by the above discussion. The proof of part (ii) is now finished.

We now give the proof of part (i). Again we consider two cases. 
First assume that $\ell_{\p_1,\p_2}\in Pa$ and $\p_2\notin \p_1^\perp$.
Then $|N_E(\p_1)\cap E_{\ell_{\mathcal{O}_i,\p_2}}|$, $1\le i\le 2$, are odd
if $q\equiv \pm 1\pmod 8$; otherwise, they are even by Lemma~\ref{intersection} (i). 
Let $\ell_s\in Se_{\p_2}$ and $\ell_s\not=\ell_{\mathcal{O}_i,\p_1}$ for $1\le i \le 2$. 
Note that $\ell_s^\perp\notin T_1$ or $T_2$. Hence $\ell_{\p_1,\ell_s^\perp}$ is either 
a secant or a passant. So $|N_E(\p_1)\cap E_{\ell_s}|$ is even by Lemma~\ref{intersection} (ii).
Set $L:=Se_{\p_1}\setminus\{\ell_{\mathcal{O}_1,\p_2}, \ell_{\mathcal{O}_2,\p_2}\}$. 
Then $|L|=\frac{q-1}{2}-2$ and
\begin{displaymath}
\begin{array}{llll}
{} & |N_E(\p_1)\cap N_E(\p_2)|\\
= &
\begin{cases}
\displaystyle\sum_{\ell\in L}|N_E(\p_1)\cap(E_{\ell}\setminus\{\p_2\})|+n_1+n_2 & \text{if}\;
q\equiv 1\pmod 4\\
\displaystyle\sum_{\ell\in L}|N_E(\p_1)\cap E_{\ell}|+n_1+n_2 &
\text{if}\;q\equiv 3\pmod 4\\
\end{cases}\\
\end{array}
\end{displaymath}
where $n_i=|N_E(\p_1)\cap E_{\ell_{\mathcal{O}_i,\p_2}}|$, $i=1,2$, are odd if 
$q\equiv \pm 1\pmod 8$; otherwise, these numbers are even. Note that
$|N_E(\p_1)\cap(E_{\ell}\setminus\{\p_2\})|$ and $|N_E(\p_1)\cap E_{\ell}|$, $\ell\in L$,
are all even numbers by the above discussion and the fact that
$|N_E(\p_1)\cap(E_{\ell}\setminus\{\p_2\})|=|N_E(\p_1)\cap E_{\ell}|$ as
$\ell_{\p_1,\p_2}\in Pa$. Thus $$|N_E(\p_1)\cap N_E(\p_2)|\equiv 0 \pmod 2.$$

Next we consider $\ell_{\p_1,\p_2}\in Pa$ and $\p_2\in \p_1^\perp$. 
This case can happen only when $q\equiv 3\pmod 4$ by applying 
Lemma~\ref{meet} (iii) to the incidence pair $(\p_1,\ell_{\p_1,\p_2})$. 
If $\ell_s\in Se_{\p_2}$ and $\ell_s\not=\p_1^\perp$, then
$\ell_{\p_1,\ell_s^\perp}$ is either a secant or a passant. Thus
$|N_E(\p_1)\cap E_{\ell_s}|$ is even by Lemma~\ref{intersection} (ii).
If $\ell_s=\p_1^\perp$, then $|N_E(\p_1)\cap E_{\ell_s}| = 0$ by part 
(iv) of Lemma~\ref{meet}. Therefore, $|N_E(\p_1)\cap N_E(\p_2)|\equiv 0\pmod 2$.

Part (iii) can be proved in the same fashion. We omit the details. \QED

\subsection{The Proof of Theorem~\ref{matrix}}

We are now ready to prove Theorem~\ref{matrix} by using the
geometric results obtained in the previous two subsections.

\vspace{0.1in}

\noindent{\bf Proof of Theorem~\ref{matrix}.}  Note that the row of 
$\B^\top \B=\B^2\pmod 2$ indexed by the point $\p_i$ can be viewed as the 
characteristic vector of $N_E(\p_i)$ and the row of $\B^3\pmod 2$ 
indexed by $\p_i$ is given by
$$\left(|N_E(\p_i)\cap E_{\ell}|\pmod 2\right)_{\ell\in Se}.$$
By Lemma~\ref{intersection}, we see that $|N_E(\p_i)\cap E_{\ell}|$ 
is odd if and only if either $(1)$ $\p_i\in \ell$ and $\ell\in Se$ for all 
$q$ or $(2)$ $\p_i\not\in \ell$, $\p_i\in T_1$ or $T_2$ (where $T_1$ 
and $T_2$ are the two tangents through $\ell^{\perp}$) and 
$q\equiv \pm 1 \pmod 8$. Therefore $\B^3\equiv \B \pmod 2$ if
$q\equiv \pm 3 \pmod 8$.

Furthermore, the row of $\B^4=\B^2\B^2\pmod 2$ indexed by $\p_i$ 
is given by $$(|N_E(\p_i)\cap N_E(\p_j)|\pmod 2)_{\p_j\in E}.$$
When $\p_i\not=\p_j$, by Lemma~\ref{intersection2} we know that 
if $q\equiv\pm 1\pmod 8$ then
\begin{equation}\label{inset1}
|N_E(\p_i)\cap N_E(\p_j)|=
\begin{cases}
d_{ij}, & \text{if}\;\ell_{\p_i,\p_j}\in Pa,\\
d_{ij}^{'}, & \text{if}\;\ell_{\p_i,\p_j}\in Se\;\text{or}\;T,
\end{cases}
\end{equation}
where $d_{ij}$ is even and $d_{ij}^{'}$ is odd;
if $q\equiv \pm3 \pmod 8$ then
\begin{equation}\label{inset11}
|N_E(\p_i)\cap N_E(\p_j)|=
\begin{cases}
b_{ij}, & \text{if}\;\ell_{\p_i,\p_j}\in Pa\;\text{or}\;T,\\
b_{ij}^{'}, & \text{if}\;\ell_{\p_i,\p_j}\in Se,
\end{cases}
\end{equation}
where $b_{ij}$ is even and $b_{ij}^{'}$ is odd.
Also, it is clear that
\begin{equation}\label{inset2}
|N_E(\p_j)| =
\begin{cases}
\frac{(q-1)(q-3)}{4}, & \text{if}\;q\equiv 1\pmod 4,\\
\frac{(q-1)(q-3)}{4} + 1, & \text{if}\;q\equiv 3 \pmod 4.\\
\end{cases}
\end{equation}
From (\ref{inset1}) and (\ref{inset2}), we see that if $q\equiv \pm
1\pmod 8$ then the row of $\B^4\pmod 2$ indexed by $\p_i$ can be
viewed as the sum of the $\Ff_2$-characteristic vector of
$N_E(\p_i)$ and that of $T_E(\p_i)$, where $T_E(\p_i):=\{\q\in E\mid
\q\in \ell,\;\ell\in T_{\p_i}\}\setminus\{\p_i\}$. And from
(\ref{inset11}) and (\ref{inset2}) we see that if $q\equiv \pm
3\pmod 8$ then the row of $\B^4\pmod 2$ indexed by $\p_i$ can be
viewed as the characteristic vector of $N_E(\p_i)$. Therefore, if
$q\equiv \pm 1\pmod 8$ then the row of $\B^5=\B^4\B\pmod 2$ indexed
by $\p_i$ is given by
$$\left(|(N_E(\p_i)\cup T_E(\p_i))\cap E_{\ell}|\pmod 2\right)_{\ell\in Se};$$
if $q\equiv \pm 3\pmod 8$ then the row of $\B^5=\B^4\B\pmod 2$ indexed by
$\p_i$ is given by $$\left(|N_E(\p_i)\cap E_{\ell}|\pmod 2\right)_{\ell\in Se},$$
where
\begin{displaymath}
|N_E(\p_i)\cap E_{\ell}|\equiv
\begin{cases}
1 \pmod 2, &\text{if}\; \p_i\in\ell,\\
0 \pmod 2, &\text{if}\; \p_i\notin \ell,\\
\end{cases}
\end{displaymath}
by the discussion in the first paragraph.

Assume that $q\equiv \pm 1\pmod 8$. Let $\mathcal{O}_1$ and
$\mathcal{O}_2$ be two points of $\mathcal{O}$ on $\p_i^\perp$.
Since $N_E(\p_i)\cap E_{\ell}$ and $T_E(\p_i)\cap E_{\ell}$ are
disjoint, by Lemma~\ref{intersection}, we have
\begin{displaymath}
\begin{array}{llll}
{} & |(N_E(\p_i)\cup T_E(\p_i))\cap E_{\ell}| \\
= & |N_E(\p_i)\cap E_{\ell}|
 + |T_E(\p_i)\cap E_{\ell}|\\
= & \left\{\begin{array}{cllllllllll}
\frac{q-3}{2} & \text{if}\; \p_i\in \ell\;\text{and}\;q\equiv 1\pmod 8& (a_1)\\
\frac{q-1}{2} & \text{if}\; \p_i\in \ell\;\text{and}\;q\equiv -1\pmod 8& (a_2)\\
|N_E(\p_i)\cap E_{\ell}|+2 & \text{if}\; \p_i\notin \ell,\;\mathcal{O}_1\notin
\ell,\;\mathcal{O}_2\notin\ell& (a_3)\\
|N_E(\p_i)\cap E_{\ell}|+1 & \text{if}\;\p_i\notin \ell,\;\mathcal{O}_1\in\ell,
\;\mathcal{O}_2\notin\ell & (a_4)\\
|N_E(\p_i)\cap E_{\ell}|+1 & \text{if}\;\p_i\notin \ell,\;\mathcal{O}_1\notin
\ell,\;\mathcal{O}_2\in\ell & (a_5)\\
|N_E(\p_i)\cap E_{\ell}| & \text{if}\;\p_i\notin \ell,\;\mathcal{O}_1\in \ell,
\;\mathcal{O}_2\in\ell & (a_6)\\
\end{array}\right.\\
\equiv &
\begin{cases}
1\pmod 2 & \text{if}\;\p_i\in \ell\\
0 \pmod 2 & \text{if}\;\p_i\notin \ell.\\
\end{cases}
\end{array}
\end{displaymath}
In Case $(a_3)$, $|N_E(\p_i)\cap E_{\ell}|$ is even by (ii) of 
Lemma~\ref{intersection} since $\p_i$ cannot be on any of the two 
tangent lines through $\ell^\perp$; in Case $(a_4)$ (respectively, 
Case $(a_5)$), we have $|N_E(\p_i)\cap E_{\ell}|$ is odd by
Lemma~\ref{intersection} (i) since $\p_1$ is on the tangent line 
$\mathcal{O}_1^\perp$ (respectively, $\mathcal{O}_2^\perp$) 
through $\ell^\perp$; in Case $(a_6)$, $|N_E(\p_i)\cap E_{\ell}|=\frac{q-1}{2}$ 
or $0$, which is even by Lemma~\ref{meet} since $\ell=\p_i^\perp$.

Therefore, $\B^5\equiv \B \pmod 2$. The proof is complete.
\QED
\begin{Remark}
In the rest of the paper, given $\p\in E$, we will use $T_E(\p)$ to denote
the set of external points excluding $\p$ that are on the tangent lines
through $\p$.
\end{Remark}

\begin{Definition} Let $\p\in E$. We define
\begin{displaymath}
{N_E(\p)}^{\rm a}=
\begin{cases}
N_E(\p)\cup \{\p\}\cup T_E(\p), & \text{if}\;q\equiv 1\pmod 8,\\
N_E(\p)\cup\{\p\}, & \text{if}\;q\equiv 5\pmod 8,\\
(N_E(\p)\cup T_E(\p))\setminus \{\p\}, & \text{if}\;q\equiv 7\pmod 8,\\
N_E(\p)\setminus\{\p\},& \text{if}\;q\equiv 3\pmod 8.\\
\end{cases}
\end{displaymath}
\end{Definition}
\begin{Corollary}\label{nullspace}
Let $F$ be an algebraic closure of $\Ff_2$. Viewing $\B$ as a matrix
with entries in $F$, we have that the characteristic vectors of
${N_{E}(\p)}^{\rm a}$ with $\p\in E$ span the null space of $\B$
over $F$.
\end{Corollary}

{\Proof} First we prove that the $F$-null space of $\B$ is equal to
the span of the rows of $\B^4+I$. If $\bf{x}$ is in the $F$-null
space of $\B$, then ${\bf x}(\B^4+ I) = \bf{x}$. That is, $\bf{x}$
is in the $F$-span of the rows of $\B^4+I$. On the other hand, if
$\bf{x}$ is in the $F$-span of the rows of $\B^4+I$, then ${\bf
y}(\B^4+I)=\bf{x}$ for some $\bf{y}$. Therefore, ${\bf
y}(\B^4+I)\B=\bf{x}\B=\bf{y}(\B^5+\B)=\bf{0}$. That is, $\bf{x}$ is
in the $F$-null space of $\B$. From the proof of
Theorem~\ref{matrix}, it is easily seen that the rows of $\B^4+I$
can be realized as the characteristic vectors of ${N_E(\p_i)}^{\rm
a}$ for $\p_i\in E$. Therefore the corollary follows. \QED

\section{The Conjugacy Classes of $H$ and Related Intersection Properties}
In this section, we give detailed information about the conjugacy classes of
$H$ and discuss their intersections with some special subsets of $H$.
\subsection{Conjugacy classes}
Recall that \begin{displaymath} H=
\left.\left\{\left(\begin{array}{ccc}
a^2 & ab & b^2 \\
2ac & ad+bc & 2bd \\
c^2 & cd & d^2 \\
\end{array}\right)\right| a, b, c, d\in \Ff_q,\;ad-bc=1\right\}
\end{displaymath}
is the subgroup of $G$ that is isomorphic to $PSL(2,q)$. Therefore, the 
conjugacy classes of $H$ can be deduced from those of $PSL(2,q)$. 
The conjugacy classes of $PSL(2,q)$ are given in the following in
terms of $2\times 2$ matrices. We refer the reader to \cite{jordan}
or \cite{schur} for the detailed calculations.
\begin{Lemma}\label{cpsl}$($\cite{jordan}, \cite{schur}$)$
The conjugacy classes of $PSL(2,q)$ can be explained
in terms of $2\times2$ matrices as follows:
If $q\equiv 1\pmod 4$ $($respectively, $q\equiv 3\pmod 4$$)$, then $D$, $[0]$, $F^+$, 
$F^-$, $[\theta_i]$ with $1\le i\le \frac{q-5}{4}$ $($respectively, $1\le i\le \frac{q-3}{4}$$)$,
$[\pi_k]$ with $1\le k\le \frac{q-1}{4}$ $($respectively, $1\le k\le \frac{q-3}{4}$$)$ are all 
conjugacy classes of $PSL(2,q)$, where $D=\left\{\pm\left(\begin{smallmatrix}1&0\\0&1
\end{smallmatrix}\right)\right\}$, $[0]$ is the class whose representative is 
$\pm\left(\begin{smallmatrix}0& -1 \\ 1& 0\end{smallmatrix}\right)$,
$F^+$ is the class whose representative is
$\pm\left(\begin{smallmatrix}1& 0 \\ 1& 1\end{smallmatrix}\right)$,
$F^-$ is the class whose representative is
$\pm\left(\begin{smallmatrix}1& 0 \\ \xi& 1\end{smallmatrix}\right)$,
$[\theta_i]$ is the class whose representative is $\pm\left(\begin{smallmatrix}t_i& 0\\ 0& 
{t_i^{-1}}\end{smallmatrix}\right)$
for some $t_i\in \Ff_q^*\setminus\{\pm 1\}$ such that $0\not=(t_i+{t_i^{-1}})^2=\theta_i$,
and $[\pi_k]$ is the class whose representative is 
$\pm\left(\begin{smallmatrix}t_k& -1\\ 1& 0 \end{smallmatrix}\right)$
for some $t_k\in \Ff_q$ such that $\pi_k=t_k^2$ and
$\pi_k-4\in\Nsq_q$.
\end{Lemma}

Let
\begin{equation}\label{defg}
g=\left(\begin{array}{ccc}
a^2 & ab & b^2 \\
2ac & ad+bc & 2bd \\
c^2 & cd & d^2 \\
\end{array}\right) \in H.
\end{equation}
For convenience, we define a map $T$ from $H$ to $\Ff_q$ by setting $T(g)=\text{tr}(g)+1$,
where $\text{tr}(g)$ is the trace of $g$;
explicitly, we have $T(g)=(a+d)^2$. If we still use $D$, $[0]$, $F^+$, $F^-$, $[\theta_i]$ 
and $[\pi_k]$ to denote the sets of images of the elements in
the corresponding classes of $PSL(2,q)$ under the isomorphism $\tau$ discussed at
the beginning of Section 2.2,
respectively, they form the conjugacy classes of $H$.

\begin{Lemma}\label{classes}
The conjugacy classes of $H$ are given as follows.
\begin{itemize}
\item[(i)] $D=\{{\bf d}(1,1,1)\}$;
\item[(ii)] $F^{+}$ and $F^{-}$, where $F^{+}\cup F^{-} = \{g\in H\mid T(g) = 4,\;g\not=\{{\bf d}(1,1,1)\}$;
\item[(iii)] $[\theta_i] =\{g\in H \mid T(g) = \theta_i\}$, $1\le i \le \frac{q-5}{4}$ if $q\equiv 1
\pmod 4$, or $1\le i\le \frac{q-3}{4}$ if $q\equiv 3 \pmod 4$, where
$\theta_i\in \Sq_q$, $\theta_i\not= 4$, and $\theta_i-4\in\Sq_q$;
\item[(iv)] $[0]=\{g\in H\mid T(g) =0\}$;
\item[(v)] $[\pi_k]=\{g\in H\mid T(g) = \pi_k\}$, $1\le k\le \frac{q-1}{4}$ if $q\equiv 1\pmod 4$, or
$1\le k\le \frac{q-3}{4}$ if $q\equiv 3 \pmod 4$, where $\pi_i\in\Sq_q$, $\pi_k\not=4$, and
$\pi_k-4\in\Nsq_q$.
\end{itemize}
\end{Lemma}
\begin{Remark}\label{order}
The set $F^+\cup F^-$
forms one conjugacy class of $G$, and splits into two equal-sized classes 
$F^{+}$ and $F^{-}$ of $H$. For our purpose, we denote $F^+\cup F^-$ by $[4]$. 
Also, each of $D$, $[\theta_i]$, $[0]$, and $[\pi_k]$ forms a single conjugacy 
class of $G$. The class $[0]$ consists of all the elements of order $2$ in
$H$.
\end{Remark}
In the following, for convenience, we frequently use $C$ to denote
any one of $D$, $[0]$, $[4]$, $[\theta_i]$, or
$[\pi_k]$. That is,
\begin{equation}\label{defC}
C=D, [0], [4], [\theta_i], \;{\rm or}\;[\pi_k].
\end{equation}

\subsection{Intersection properties}
We study
the intersection sizes of certain subsets of $H$ with the conjugacy classes of $H$.

\begin{Definition}\label{definition} 
Let $\p,\q\in E$ be two external points, $\ell$ a secant line,
and $W\subseteq E$. We define
\begin{displaymath}
\begin{array}{lll}
\mathcal{H}_{\p,\q} & = &\{h\in H\mid (\p^\perp)^h\in Se_{\q}\},\\
\mathcal{S}_{\p,\ell} & = &\{h\in H\mid (\p^\perp)^h=\ell\},\\
\mathcal{U}_{\p, W} & =& \{h\in H\mid \p^h\in W\}.\\
\end{array}
\end{displaymath}
That is, $\mathcal{H}_{\p,\q}$ consists of all the elements in $H$
that map the secant line $\p^\perp$ to a secant line through $\q$,
$\mathcal{S}_{\p,\ell}$ is the set of elements in $H$ that map
$\p^\perp$ to the secant line $\ell$, and $\mathcal{U}_{\p, W}$ is
the set of elements in $H$ that map $\p$ to a point in $W$.
\end{Definition}

The following lemma and corollary are clear.
\begin{Lemma}\label{group_1}
Let $g\in G$, $\p,\q\in E$, $W\subseteq E$, and $\ell$ a secant
line. Then
$\mathcal{H}_{\p,\q}^g=\mathcal{H}_{\p^g,\q^g}$,
$\mathcal{S}_{\p,\ell}^g = \mathcal{S}_{\p^g,\ell^g}$,
and
$\mathcal{U}_{\p, W}^g=\mathcal{U}_{\p^g, W^g}$,
where $\mathcal{H}_{\p,\q}^g
=\{g^{-1}h g\mid h \in \mathcal{H}_{\p,\q}\}$, $\mathcal{S}_{\p,\ell}^g=\{g^{-1} h g\mid
h\in \mathcal{S}_{\p,\q}\}$, and $\mathcal{U}_{\p, W}^g=\{g^{-1}h g\mid h \in \mathcal{U}_{\p, W}\}$.
\end{Lemma}
\begin{Corollary}\label{group_intersection}
Let $g\in G$ and $C$ be given in $($\ref{defC}$)$ and let $\p$, 
$\q\in E$, $W\subseteq E$, and $\ell$ a secant line. Then 
$(C\cap\mathcal{H}_{\p,\q})^g =C\cap\mathcal{H}_{\p^g, \q^g}$, 
$(C\cap\mathcal{S}_{\p,\ell})^g=C\cap \mathcal{S}_{\p^g,\ell^g}$, 
and
$(\C\cap\mathcal{U}_{\p,W})^g=C\cap\mathcal{U}_{\p^g,W^g}$.
\end{Corollary}

In the following lemmas, we investigate the parity of
$|\mathcal{H}_{\p,\q}\cap C|$ for any two external points $\p$ and
$\q$. This information will be used in the proof of Lemma~\ref{y4}.

\begin{Lemma}\label{stabilizer_intersection}
Let $\p\in E$ and $K=Stab_H(\p)$. Then
\begin{itemize}
\item[(i)] $|K\cap D| = 1$.
\item[(ii)] $|K\cap [0]|=\frac{q+1}{2}$ or $\frac{q-1}{2}$ according as
$q\equiv 1$ or $3\pmod 4$.
\item[(iii)] $|K\cap [\pi_k]| = 0$ for each $k$.
\item[(iv)] $|K\cap [\theta_i]|=2$ for each $i$.
\item[(v)] $|K\cap [4]| =0$.
\end{itemize}
\end{Lemma}
{\Proof} Part (i) is obvious. Let $\q$ be an arbitrary external
point. By Lemma~\ref{transitive}, there exists an element 
$g\in G$ such that $\p=\q^g$. Then 
$Stab_H(\p)=Stab_H(\q^g)=Stab_H(\q)^g$ 
by Lemma~\ref{stabilizer}. Therefore 
$$|Stab_H(\p)\cap C|=|(Stab_H(\q)\cap C)^g|=|Stab_H(\q)\cap C|,$$ 
where $C$ is given in (\ref{defC}). It follows that we may 
assume that $\p=(0,1,0)$. Let $g\in H$ be given by (\ref{defg}).
Then the element $g$ given by $(a,b,c,d)$ is in $K\cap C$ if and only
if the following system of equations holds.
\begin{displaymath}
\begin{array}{cccccc}
ac & = & 0 \\
ad+bc & = & u \\
bd & = & 0\\
ad-bc & = & 1 \\
a+d & = & s\\
\end{array}
\end{displaymath}
where $u\in \Ff_q^*$, and $s=0$, $\pm 2$, $\pm\sqrt{\theta_i}$, or
$\pm\sqrt{\pi_k}$. Solving this system of equations, we have either
$a=d=0$ or $b=c=0$.

If $a=d=0$, then $bc=-1$ and $s=0$. In this case, we obtain
$\frac{q-1}{2}$ elements $g$ in $K\cap [0]$.

If $b=c=0$, then $g={\bf d}(a^2,1,a^{-2})$, so either $g={\bf d}(1,1,1)$ or 
$g \in [\theta_i]$, where $\theta_i=a^2+a^{-2}+2$, or $q\equiv 1 \pmod 4$
and $g={\bf d}(-1,1,-1)\in [0]$. Since $a$ and $a^{-1}$ give the same
$\theta_i$, there are exactly $2$ elements $g$ for each $\theta_i$ .
\QED

Recall that $\ell_{\p,\q}$ denotes the line through the points $\p$
and $\q$, and $T_{\p}$ denotes the set of tangent lines through
$\p$.
\begin{Lemma}\label{line_1}
Assume that $q\equiv 1\pmod 4$. Let $\p$ and $\q$ be two distinct
external points and let $C=D$, $[4]$, $[\pi_k]$ $($$1\le k\le
\frac{q-1}{4}$$)$, or $[\theta_i]$ $($$1\le i\le \frac{q-5}{4}$$)$.
\begin{itemize}
\item[(i)] Suppose that $\ell_{\p,\q}\in Se_{\p}$ and $\q\notin \p^\perp$. 
If $|\mathcal{H}_{\p,\q}\cap C|$ is odd, then $C$ must be equal to
one of two distinct classes of type $[\theta_{i}]$.
\item[(ii)] Suppose that $\ell_{\p,\q}\in Se_{\p}$ and $\q\in \p^\perp$. 
If $|\mathcal{H}_{\p,\q}\cap C|$ is odd, then $C=D$.
\item[(iii)] Suppose that $\ell_{\p,\q}\in Pa_{\p}$. Then 
$|\mathcal{H}_{\p,\q}\cap C|$ is always even.
\item[(iv)] Suppose that $\ell_{\p,\q}\in T_{\p}$. If 
$|\mathcal{H}_{\p,\q}\cap C|$ is odd, then 
$C=[\theta_i]$ fomr some $1\le i \le \frac{q-5}{4}$ $($possibly
more than one $[\theta_i]$ intersect with $\mathcal{H}_{\p,\q}$
in an odd number of group elements$)$.
\end{itemize}
\end{Lemma}

{\Proof} Without loss of generality we may assume that
$\p=(0,1,0)$ since $G$ acts transitively on $E$. Let
$K=Stab_G(\p)$. Assume that $\ell_1$ and $\ell_2$ 
are two lines in $Pa_{\p}$, $Se_{\p}$, or $T_{\p}$, and 
$\q\in \ell_1$. 
By Proposition~\ref{K_transitive}, it follows that
$\ell_1^g=\ell_2$ for some $g\in K$, and so $\q^g\in \ell_2$.
Moreover, Corollary~\ref{group_intersection} gives
$$|\mathcal{H}_{\p,\q}\cap C|=|(\mathcal{H}_{\p,\q}\cap C)^g|
= |\mathcal{H}_{\p^g,\q^g}\cap C|= |\mathcal{H}_{\p,\q^g}\cap C|.$$
Therefore, to prove the lemma, it is enough to take $\q$ to be an
arbitrary external point on a {\it special} secant or passant or
tangent line through $\p$.

\vspace{0.1in}

\noindent{(i)} $\q\notin \p^\perp=[0,1,0]$ and 
$\ell_{\p,\q}=[1,0,y]\in Se_{\p}$, for some $y\in\Sq_q$.

In this case, we have $\q=(1,m,-y^{-1})$ for some
$m\not=0$ and $m^2+y^{-1}\in\Sq_q$. From the computations 
done in the proof of Lemma~\ref{intersection}, we know that 
both $\ell_1=[1,0,y]$ and $\ell_2=[1,-2m^{-1},-y]$ are fixed by
$Stab_K(\q)$, and $Se_{\q}\setminus\{\ell_1, \ell_2\}$ splits into
$\frac{q-5}{4}$ orbits of length $2$ under the action of
$Stab_K(\q)$. Let $\mathcal{R}$ be a set of the representatives of
these orbits of length $2$. Then by
Corollary~\ref{group_intersection}, we have
\begin{displaymath}
\begin{array}{lll}
|\mathcal{H}_{\p, \q}\cap C| & = & \displaystyle\sum_{\ell\in
Se_{\q}}|\mathcal{S}_{\p,\ell}
\cap C| \\
{} & = & |\mathcal{S}_{\p,\ell_1}\cap C| +
|\mathcal{S}_{\p,\ell_2}\cap C| +
\displaystyle\sum_{\ell\in \mathcal{R} }2|\mathcal{S}_{\p,\ell}\cap C|.\\
\end{array}
\end{displaymath}
Here we have used the fact that if $\{\ell,\ell'\}$ is an orbit of
secant lines through $\q$, then $|\mathcal{S}_{\p,\ell}\cap
C|=|\mathcal{S}_{\p,\ell'}\cap C|$. From the above equation we see
that the parity of $|\mathcal{H}_{\p,\q}\cap C|$ is the same as that
of $|\mathcal{S}_{\p,\ell_1}\cap C| + |\mathcal{S}_{\p,\ell_2}\cap
C|$. In the following, we determine the parity of
$|\mathcal{S}_{\p,\ell_1}\cap C|$ and $|\mathcal{S}_{\p,\ell_2}\cap
C|$.

$(1a)$ Let $g\in \mathcal{S}_{\p,\ell_1}\cap C$, where $g$ is given
by (\ref{defg}). Since $(\p^\perp)^g$ is determined by the the following
column vector
\begin{displaymath}
\left(\begin{array}{ccc}
d^2 & -bd & b^2 \\
-2cd & ad+bc & -2ab \\
c^2 & -ac & a^2 \\
\end{array}
\right)
\left(
\begin{array}{ccc}
0 \\ 1 \\ 0
\end{array}
\right)=
\left(
\begin{array}{c}
-bd \\ ad+bc \\ -ac
\end{array}
\right),
\end{displaymath}
we see that the quadruple $(a,b,c,d)$ determining $g$ satisfies the
following system of equations
\begin{equation}\label{eq1}
\begin{array}{ccccc}
ad+bc & = & 0\\
ac& = & ybd \\
ad-bc & = & 1 \\
a+d & = & s,
\end{array}
\end{equation}
where $s=\pm 2,\pm \sqrt{\pi_k},\pm\sqrt{\theta_i}$. Therefore
$|\mathcal{S}_{\p,\ell_1}\cap C|$ can be determined by the number of
solutions $(a,b,c,d)$ to the equations in (\ref{eq1}). The equations
in (\ref{eq1}) yield $a^2-sa+\frac{1}{2}=0$, whose discriminant is
$s^2-2$. If $s^2-2=0$ (so $2$ is a square), then $a=d=\frac{s}{2}$
and $b=c=\pm\sqrt{-\frac{1}{2y}}$, giving two distinct group
elements $g\in [2]$. If $s^2-2\in\Sq_q$, then the equations in
$(\ref{eq1})$ yield $0$ or $4$ different solutions of $(a,b,c,d)$ such
that $(a,b,c,d)$ and $(-a,-b,-c,-d)$ appear at the same time if any; 
these solutions give rise to $2$ distinct group elements $g$ in $[s^2]$.
If $s^2-2\in\Nsq_q$, it is obvious that
$|\mathcal{H}_{\p,\ell_1}\cap [s^2]| =0$. It is also clear that
$|\mathcal{S}_{\p,\ell_1}\cap D| =0$. Therefore,
$|\mathcal{S}_{\p,\ell_1}\cap C|$ is even for all choices of $C$ as
specified in the statement of the lemma.

$(1b)$ Let $g\in \mathcal{S}_{\p,\ell_2}\cap C$, where $g$ is given
by (\ref{defg}). Similarly, the quadruple $(a,b,c,d)$ determining
$g$ satisfies the following system of equations
\begin{equation}\label{eqq2}
\begin{array}{cccccccc}
ad+bc & = & \frac{2}{m}(bd) \\
ac & =  & -y(bd)   \\
ad-bc & = & 1 \\
a+d & = & s,
\end{array}
\end{equation}
where $s=\pm 2,\pm\sqrt{\pi_k},\pm\sqrt{\theta_i}$. From the first
two equations in $(\ref{eqq2})$ we obtain
$$(ad-bc)^2=4(\frac{1}{m^2}+y)(bd)^2;$$ that is,
$bd=\pm\frac{1}{\sqrt{w}}$, where $w=4(\frac{1}{m^2}+y)\in
\Sq_q\;\;\text{and}\;\; y\in\Sq_q$. Thus,
$d^2-ds+(\frac{1}{m\sqrt{w}}+\frac{1}{2}) = 0$ or
$d^2-ds+(-\frac{1}{m\sqrt{w}}+\frac{1}{2}) = 0$. The discriminants
of these two quadratic equations are
$$\Delta_1(s^2)=s^2-4(\frac{1}{m\sqrt{w}}+\frac{1}{2})$$ and
$$\Delta_2(s^2) =s^2-4(-\frac{1}{m\sqrt{w}}+\frac{1}{2}),$$ respectively.
Note that neither $4(\frac{1}{m\sqrt{w}}+\frac{1}{2})$ nor
$4(-\frac{1}{m\sqrt{w}}+\frac{1}{2})$ can be $0$ since 
$y\not=0$ and
\begin{equation}\label{chill}
(4)(\frac{1}{m\sqrt{w}}+\frac{1}{2})(4)(-\frac{1}{m\sqrt{w}}+\frac{1}{2})
=\frac{16y}{w}\in \Sq_q.
\end{equation}

If $|\mathcal{S}_{\p, \ell_2}\cap [s^2]|$ is odd, then either
$\Delta_1(s^2)=0$ or $\Delta_2(s^2)=0$ since
$(\pm\frac{1}{m\sqrt{w}}+\frac{1}{2})\not=0$. It follows that either
$s^2=4(\frac{1}{m\sqrt{w}}+\frac{1}{2})$ or
$s^2=4(-\frac{1}{m\sqrt{w}}+\frac{1}{2})$, and so by (\ref{chill}), we
see that $s^2$ must be any one of  $\theta_{i_1}$ and $\theta_{i_2}$, 
where $\theta_{i_1}=4(\frac{1}{m\sqrt{w}}+\frac{1}{2})$,
$\theta_{i_2}=4(-\frac{1}{m\sqrt{w}}+\frac{1}{2})$.

Combing $(1a)$ and $(1b)$, we see that, if $|\mathcal{H}_{\p,\q}\cap
C|$ is odd, then $C$ must be any one of two classes $[\theta_{i_1}]$ and
 $[\theta_{i_2}]$, where $\theta_{i_1}, \theta_{i_2}$ are given as above.

\vspace{0.1in} \noindent{(ii)}
$\q=(1,0,-y^{-1})\in \p^\perp$ and $\ell_{\p,\q}=[1,0,y]\in Se_{\p}$, 
for some $y\in\Sq_q$.

In this case, we know that both $\p^\perp$ and
$\ell_{\p,\q}=[1,0,y]$ are fixed by $Stab_{K}(\q)$, and 
$Se_{\q}\setminus\{\ell_{\p,\q},\p^\perp\}$ splits into 
$\frac{q-5}{4}$ orbits of length $2$ under the action of
$Stab_K(\q)$ from the computations done in the proof of
Lemma~\ref{intersection}. Let $\mathcal{R}$ be a set of the
representatives of these orbits of length $2$. Then we have
\begin{displaymath}
\begin{array}{lll}
|\mathcal{H}_{\p,\q}\cap C| & = &
\displaystyle\sum_{\ell\in Se_{\q}}|\mathcal{S}_{\p,\ell}\cap C| \\
{} & = & |\mathcal{S}_{\p,\ell_{\p,\q}}\cap C| +|K\cap C| +
\displaystyle\sum_{\ell\in \mathcal{R}}2|\mathcal{S}_{\p, \ell}\cap C|.\\
\end{array}
\end{displaymath}
Here we have used the fact that
$\mathcal{S}_{\p,\p^{\perp}}=Stab_G(\p^{\perp})=K$ by
Corollary~\ref{commutative}. From this equation we see that
$|\mathcal{H}_{\p,\q}\cap C|$ is even except when $C=D$ since
$|\mathcal{S}_{\p,\ell_{\p,\q}}\cap C| $ is even by the results in
$(1a)$ and the fact that $|K\cap C|$ is even except when $C= D$ by
Lemma~\ref{stabilizer_intersection}.

\vspace{0.1in} \noindent{(iii)} $\q\notin\p^\perp$ and
$\ell_{\p,\q}=[1,0,x]\in Pa_{\p}$, for some $x\in\Nsq_q$.

In this case, we know that $Stab_K(\q)$ has $\frac{q-1}{4}$ orbits
of length $2$ on $Se_{\q}$ by computations similar to those done in
the proof of Lemma~\ref{intersection}. Let $\mathcal{R}$ be a set of
representatives of these orbits. Then
$$|\mathcal{H}_{\p,\q}\cap C|=\displaystyle\sum_{\ell\in \mathcal{R}}2|
\mathcal{S}_{\p,\ell}\cap C|,$$
which is always even.

\vspace{0.1in} \noindent{(iv)} $\ell_{\p,\q}\in T_{\p}$.

Without loss of generality, we may take
$\ell_{\p,\q}=[1,0,0]$ and $\q=(0,1,f)$ with $f\in \Ff_q^*$. Then
  $$Se_{\q}=\{[0,1,-f^{-1}]\}\cup\{[1,-fu,u]\mid u\in \Ff_q^*,
  f^2u^2-4u\in \Sq_q\}.$$ Let $U_f=\{u\mid u\in\Ff_q^*, f^2u^2-4u\in\Sq_q\}$.
  For convenience, set $\ell_u:=[1,-fu,u]$ for $u\in U_f$
  and $\ell_3:=[0,1,-f^{-1}]$. Then $Se_{\q}=\{\ell_3\}\cup\{\ell_u\mid u\in U_f\}$.
  Using (\ref{stab}), we see easily that $Stab_{K}(\q)$ contains the identity element only.
  Therefore, to find the parity of $|\mathcal{H}_{\p, \q}\cap C|$, we need to determine
  the parity of each term in the sum
  $$|\mathcal{S}_{\p, \ell_3}\cap C|+\displaystyle\sum_{u\in U_f}|
  \mathcal{S}_{\p, \ell_u}\cap C|=|\mathcal{H}_{\p,\q}\cap C|.$$
  It is clear that $|\mathcal{H}_{\p, \q}\cap D|=0$ since $\p^\perp\notin Se_{\q}$.
  The rest is devoted to the cases where $C$ is neither $D$ nor $[0]$.

$(4a)$ Let $g\in \mathcal{S}_{\p,\ell_3}\cap C$, where $C\not= D,
[0]$, and $g$ is given by (\ref{defg}). The quadruple $(a,b,c,d)$
determining $g$ in $\mathcal{S}_{\p,\ell_3}\cap C$ satisfy the
following system of equations
  \begin{equation}\label{eq3}
  \begin{array}{cccccc}
  bd & = & 0\\
  ad+bc & = & t \\
  ac & = & \frac{t}{f} \\
  ad-bc & = & 1\\
  a+d & = & s\\
  \end{array}
  \end{equation}
where $s=\pm 2,\pm \sqrt{\theta_i},\pm \sqrt{\pi_k}$, and $t\in\Ff_q^*$.

If $d=0$, then $s\not=0$, $t=-1$, $b=fs$, $c=-\frac{1}{fs}$, and
$a=s$. In this case, we obtain a unique group element $g$ in each
$[s^2]$. If $d\not=0$, then $t=1$ and $a^2-sa+1=0$, whose
discriminant is $s^2-4$. This quadratic solution has a single
solution $a\in \Ff_q$ if and only if $s^2-4=0$. Therefore,
$|\mathcal{S}_{\p,\ell_3}\cap C|$ is odd except when $C=[4]$.

$(4b)$ Let $g\in \mathcal{S}_{\p,\ell_3}\cap C$, where $C\not= D,
[0]$, and $g$ is given by (\ref{defg}). The quadruple $(a,b,c,d)$
determining $g$ satisfy the following system of equations
\begin{equation}\label{b3}
\begin{array}{ccccc}
ad+bc & =& fubd \\
ac & = & ubd \\
ad-bc & = & 1 \\
a+d & = & s,
\end{array}
\end{equation}
where $s=\pm 2,\pm\sqrt{\pi_k},\pm\sqrt{\theta_i}$, and $u\in U_f$.
The above equations in (\ref{b3}) yield
$d^2-sd+\frac{1}{2}(\frac{fu}{r}+1)=0$
or
$d^2-sd+\frac{1}{2}(-\frac{fu}{r}+1)=0$,
where
\begin{equation}\label{rr}
r=\sqrt{f^2u^2-4u}.
\end{equation}
The discriminants of the above two quadratic equations are
$$\Delta_1(s^2,u):=s^2-2(\frac{fu}{r}+1)$$
and $$\Delta_2(s^2,u):=s^2-2(-\frac{fu}{r}+1),$$ respectively.
Note that
\begin{equation}\label{u3}
(2)(\frac{fu}{r} +1)(2)(-\frac{fu}{r}+1)= -\frac{16u}{f^2u^2-4u}
\end{equation}
 and $2(-\frac{fu}{r} +1)-4=-2(\frac{fu}{r}+1)$.

Set $U_f^+:=\{u\in \Sq_q\mid f^2u^2-4u\in \Sq_q\}$ and
$U_f^-:=\{u\in \Nsq_q\mid f^2u^2-4u \in \Sq_q\}$. Then
$U_f=U_f^+\cup U_f^-$. Moreover, it is easy to see that $|U_f^-|$ is
equal to the number of $w\in\Nsq_q$ satisfying $w-4\in\Nsq_q$ and
this number is $|(\Nsq_q-1)\cap \Nsq_q|=\frac{q-1}{4}$ 
by (i) of Lemma~\ref{cs}.
Given a line $\ell_u$ with $u\in U_f$, $|\mathcal{S}_{\p,\ell_u}\cap
[s^2]|$ is odd if and only if either $\Delta_1(s^2,u)=0$ or
$\Delta_2(s^2,u)=0$ since $\frac{1}{2}(\pm\frac{fu}{r}+1)\not=0$;
that is, $|\mathcal{S}_{\p,\ell_u}\cap [s^2]|$ is odd if and only if
either $2(\frac{fu}{r} +1)\in\Sq_q$ or $2(-\frac{fu}{r}
+1)\in\Sq_q$.

When $u\in U_f^+$ and $|\mathcal{S}_{\p,\ell_u}\cap [s^2]|$ is odd,
$s^2$ is equal to either $\theta_{i_1}$ or $\theta_{i_2}$ by
(\ref{u3}), where $\theta_{i_1}=2(\frac{fu}{r}+1)$ and
$\theta_{i_2}=2(-\frac{fu}{r}+1)$. Given $0\not=u_i\in U_f^+$,
$r_i=\sqrt{f^2u_i^2-4u_i}$, $i=1$ or $2$, we have
$\frac{fu_1}{r_1}=\pm\frac{fu_2}{r_2}$ if and only if $u_1=u_2$.
Therefore the two conjugacy classes determined by $u_1$ are
different from those determined by $u_2$ if $u_1\not=u_2$.

When $u\in U_f^-$, exactly one of $2(\frac{fu}{r}+1)$ and
$2(-\frac{fu}{r}+1)$ is a square by (\ref{u3}), thus at most one of
$\Delta_1(s^2,u)$ and $\Delta_2(s^2,u)$ can be zero. This shows that
for each $u\in U_f^-$, there is a unique class $[\pi_k]$,
$|[\pi_k]\cap \mathcal{S}_{\p,\ell_u}|$ is odd, where $\pi_k$ is one
of $2(-\frac{fu}{r} +1)$ and $2(\frac{fu}{r} +1)$ depending on which
one is a square. Given $0\not=u_i\in U_f^-$,
$r_i=\sqrt{f^2u_i^2-4u_i}$, $i=1$ or $2$, we have
$\frac{fu_1}{r_1}=\pm\frac{fu_2}{r_2}$ if and only if $u_1=u_2$.
Therefore the (unique) conjugacy class determined by $u_1$ is
difference from the one determined by $u_2$. Since there are
$\frac{q-1}{4}$ classes $[\pi_k]$ and $|U_f^-|=\frac{q-1}{4}$, it
follows that, when $u$ runs through $U_f^-$ once, each $\pi_k$ with
$1\le k\le \frac{q-1}{4}$ appears in
$$\{2(\frac{fu}{r} +1)\in\Sq_q\mid u\in U_f^-\}\cup \{2(-\frac{fu}{r} +1)
\in\Sq_q\mid u \in U_f^-\}.$$
Therefore, in each class $[\pi_k]$, there are an odd number of group
elements mapping $\p^\perp$ to $\ell_u$, where $u$ determines
$\pi_k$.

Combining $(4a)$, $(4b)$ and the relation
$$|\mathcal{H}_{\p, \q}\cap [s^2]|=|\mathcal{S}_{\p, \ell_3}\cap
[s^2]| +\displaystyle\sum_{u\in U_f^+}|\mathcal{S}_{\p, \ell_u}\cap
[s^2]| +\displaystyle\sum_{u\in U_f^-}|\mathcal{S}_{\p, \ell_u}\cap
[s^2]|,$$ we see that, if $|\h_{\p, \q}\cap [s^2]|$ is odd, then
there are an odd number of terms whose values are odd in the right
hand side of the above equation; this can occur only when
$s^2=\theta_i$ with $1\le i\le \frac{q-5}{4}$ and there are probably
more than one $\theta_i$ satisfying this condition. Part (iv) is
now proved. \QED

\begin{Lemma}\label{y_22}
Assume that $q\equiv 3\pmod 4$. Let $\p$ and $\q$ be two distinct
external points and let $C=D$, $[4]$, $[\pi_k]$ $(1\le k\le
\frac{q-3}{4})$, or $[\theta_i]$ $(1\le i \le \frac{q-3}{4})$.
\begin{itemize}
\item[(i)] Suppose that $\ell_{\p,\q}\in Se_{\p}$ and $\q\notin \p^\perp$.
Then $|\mathcal{H}_{\p,\q}\cap C|$ is always even.
\item[(ii)] Suppose that $\ell_{\p,\q}\in Pa_{\p}$ and $\q\in \p^\perp$.
If $|\mathcal{H}_{\p,\q}\cap C|$ is odd, then $C=D$.
\item[(iii)] Suppose that $\ell_{\p,\q}\in Pa_{\p}$ and $\q\notin \p^\perp$.
If $|\mathcal{H}_{\p,\q}\cap C|$ is odd, then $C$ must be equal to one of
two distinct classes of type $[\pi_{k}]$.
\item[(iv)] Suppose that $\ell_{\p,\q}\in T_{\p}$. If $|\mathcal{H}_{\p,\q}\cap C|$
is odd, then $C=[\pi_k]$ for exactly one $k$.
\end{itemize}
\end{Lemma}
{\Proof} The proof is basically identical to that of Lemma~\ref{line_1}.
We omit the details.\QED
\begin{Lemma}\label{pisq}
Let $\p\in E$. Then $|\mathcal{H}_{\p,\p}\cap C|$ is even for $C=D$, $[4]$,
$[\pi_k]$ or $[\theta_i]$.
\end{Lemma}
{\Proof} Since $G$ is transitive on $E$, again, without loss of
generality, we may assume that $\p=(0,1,0)$. Let
$K=Stab_G(\p)$. Then $K$ is transitive on $Se_{\p}$ by
Proposition~\ref{K_transitive}. Furthermore, we have
$|\mathcal{S}_{\p,\ell}\cap C|=|\mathcal{S}_{\p,\ell_1}\cap C|$ for
$\ell, \ell_1\in Se_{\p}$ as $\ell$ and $\ell_1$ are in the same
orbit of $K$ on $Se_{\p}$. Since $K$ is transitive on $Se_{\p}$, we
have
\begin{equation}\label{same}
\begin{array}{llllll}
|\mathcal{H}_{\p,\p}\cap C| & = & \displaystyle\sum_{\ell\in Se_{\p}}
|\mathcal{S}_{\p,\ell}\cap C| & = & \frac{q-1}{2}|\mathcal{S}_{\p,\ell_1}\cap C|,
\end{array}
\end{equation}
where $\ell_1$ is a fixed line in $Se_{\p}$.

If $q\equiv 1\pmod 4$, the last term in~(\ref{same}) is always even,
and so $|\mathcal{H}_{\p,\p}\cap C|$ is even for each given $C$. If
$q\equiv 3\pmod 4$, we may take $\ell_1=[1,0,y]$ for
some $y\in \Nsq_q$. The same computations as those for the case
$(1a)$ in the proof of Lemma~\ref{line_1} show that
$|\mathcal{S}_{\p,\ell_1}\cap C|$ is even for each given $C$ as
well. The proof is complete. \QED

\begin{Definition}
Let $\p\in E$. We define $N_E^{'}(\p):= E\setminus {N_E(\p)}^{\rm
a}$. Then $N_E^{'}(\p)$ is the set of external points $($excluding
$\p$$)$ on the passant lines through $\p$ if $q\equiv 1\pmod 8$, and
it is the set of external points $($excluding $\p$$)$ on either the passant
lines or the tangent lines through $\p$ if $q\equiv 5\pmod 8$.
\end{Definition}

The following two lemmas are important in the proof of
Lemma~\ref{y4}. Since the proofs of the two lemmas are
quite similar to that of Lemma~\ref{line_1} and the computations
involved are somewhat tedious, we omit the proofs.

\begin{Lemma}\label{y_1}
Assume that $q\equiv 1\pmod 4$. Let $\p$ and $\q$ be two distinct
external points and let $C=D$, $[4]$, $[\pi_k]$ $(1\le k\le
\frac{q-1}{4})$, or $[\theta_i]$ $(1\le i \le \frac{q-5}{4})$.
\begin{itemize}
\item[(i)] Suppose that $\ell_{\p,\q}\in Pa_{\p}$. If
$|\mathcal{U}_{\p,N_E^{'}(\q)}\cap C|$ is odd, then $C=D$, or
$[\pi_k]$ for exactly one $k$.
\item[(ii)] Suppose that $\ell_{\p,\q}\in Se_{\p}$.Then
$|\mathcal{U}_{\p,N_E^{'}(\q)}\cap C|$ is always even.
\item[(iii)] Suppose that $\ell_{\p,\q}\in T_{\p}$. If $q\equiv 1\pmod 8$ and
$|\mathcal{U}_{\p,N_E^{'}(\q)}\cap C|$ is odd, then $C= [\pi_k]$ for
$1\le k\le \frac{q-1}{4}$; if $q\equiv5\pmod 8$ and 
$|\mathcal{U}_{\p,N_E^{'}(\q)}\cap C|$
is odd, then $C= D$, or $[\pi_k]$ for $1\le k\le \frac{q-1}{4}$.
\end{itemize}
\end{Lemma}

\begin{Lemma}\label{y_2}
Assume that $q\equiv 3\pmod 4$. Let $\p$ and $\q$ be two distinct
external points and let $C=D$, $[4]$, $[\pi_k]$ $(1\le k\le
\frac{q-3}{4})$, or $[\theta_i]$ $(1\le i \le \frac{q-3}{4})$.
\begin{itemize}
\item[(i)] Suppose that $\ell_{\p,\q}\in Se_{\p}$.
If $|\mathcal{U}_{\p,{N_E(\q)}^{\rm a}}\cap C|$ is odd, then $C=D$,
or $[\theta_i]$ for exactly one $i$.
\item[(ii)] Suppose that $\ell_{\p,\q}\in Pa_{\p}$.
Then $|\mathcal{U}_{\p,{N_E(\q)}^{\rm a}}\cap C|$ is always even.
\item[(iii)] Suppose that $\ell_{\p,\q}\in T_{\p}$. If $q\equiv 3 \pmod 8$ and
$|\mathcal{U}_{\p,{N_E(\q)}^{\rm a}}\cap C|$ is odd, then
$C=[\theta_i]$ for $1\le i\le \frac{q-3}{4}$; if $q\equiv 7 \pmod 8$
and $|\mathcal{U}_{\p,{N_E(\q)}^{\rm a}}\cap C|$ is odd, then $C=D$,
or $[\theta_i]$ for $1\le i\le \frac{q-3}{4}$.
\end{itemize}
\end{Lemma}

\begin{Lemma}\label{y3}
Let $\p\in E$ and let $C=D$, $[0]$, $[4]$,
$[\pi_k]$, or $[\theta_i]$.
\begin{itemize}
\item[(i)] If $q\equiv 1\pmod 4$, then $|\mathcal{U}_{\p, N_E^{'}(\p)}\cap C|$
is always even.
\item[(ii)] If $q\equiv 3\pmod 4$, then $|\mathcal{U}_{\p, {N_E(\p)}^{\rm a}}\cap C|$
is always even.
\end{itemize}
\end{Lemma}
{\Proof} We know that $Stab_G(\p)$ is transitive on $Pa_{\p}$,
$Se_{\p}$, and $T_{\p} : =\{\ell_1,\ell_2\}$, where $\ell_1$ and
$\ell_2$ are the two tangent lines through $\p$.  We first consider
the case where $q\equiv 1\pmod 4$. If $q\equiv 1\pmod 8$, then
\begin{equation*}\label{a11}
\begin{array}{llll}
|\mathcal{U}_{\p, N_E^{'}(\p)}\cap C| & = & \displaystyle\sum_{\ell\in Pa_{\p}}
|\mathcal{U}_{\p, E_{\ell}\setminus\{\p\}}\cap C| &{}\\\
{} & = & \frac{q-1}{2} |\mathcal{U}_{\p,
E_{\ell^*}\setminus\{\p\}}\cap C|,&{}
\end{array}
\end{equation*}
where $\ell^*$ is any fixed passant line through $\p$. Note that
$\frac{q-1}{2}$ is even in this case. We see that $|\mathcal{U}_{\p,
N_E^{'}(\p)}\cap C|$ is even as claimed.

If $q\equiv 5\pmod 8$, then
\begin{equation*}
\begin{array}{llll}
|\mathcal{U}_{\p, N_E^{'}(\p)}\cap C| & = &
\displaystyle\sum_{\ell\in Pa_{\p}} |\mathcal{U}_{\p,
E_{\ell}\setminus\{\p\}}\cap C| +|\mathcal{U}_{\p,
E_{\ell_1}\setminus\{\p\}}\cap C|
+|\mathcal{U}_{\p, E_{\ell_2}\setminus\{\p\}}\cap C| & {}\\\
{} & = & \frac{q-1}{2} |\mathcal{U}_{\p,
E_{\ell^*}\setminus\{\p\}}\cap C|+ 2|\mathcal{U}_{\p,
E_{\ell_1}\setminus\{\p\}}\cap C|,&{}
\end{array}
\end{equation*}
where $\ell^*$ is any fixed passant line through $\p$. It is
easily seen that $|\mathcal{U}_{\p, N_E^{'}(\p)}\cap C|$ is even.

The conclusion in the case where $q\equiv 3\pmod 4$ can be similarly
proved. We omit the details.
\QED

\section{Group Algebra $FH$}
\subsection{2-Blocks of H}
Recall that $H\cong PSL(2,q)$. In this section we recall several
results on the $2$-blocks of $H$. The results and statements in
this section are standard in the theory of characters and blocks of
finite groups. We refer the reader to {\rm \cite{gabriel}} or {\rm
{\cite{brauer}}} for a general introduction on this subject.

Let $\mathbf{R}$ be the ring of algebraic integers in the complex
field $\Cc$. We choose a maximal ideal $\mathbf{M}$ of $\mathbf{R}$
containing $2\mathbf{R}$. Let $F=\mathbf{R}/\mathbf{M}$ be the
residue field of characteristic $2$, and let $* :
\mathbf{R}\rightarrow F$ be the natural ring homomorphism. Define
\begin{equation}
\begin{array}{llll}\label{ring_s}
\mathbf{S}& =& \{\frac{r}{s}\mid r\in\mathbf{R},\;s\in\mathbf{R}\setminus\mathbf{M}\}.\\
\end{array}
\end{equation}
Then it is clear that the map $* : \mathbf{S}\rightarrow F$ defined
by $(\frac{r}{s})^* = r^*(s^*)^{-1}$ is a ring homomorphism with
kernel $\mathcal{P} = \{\frac{r}{s}\mid
r\in\mathbf{M},\;s\in\mathbf{R}\setminus\mathbf{M}\}$. In the rest
of this paper, $F$ will always be the field of characteristic 2
constructed as above. Note that $F$ is an algebraic closure of
$\Ff_2$.

Let $Irr(H)$ and $IBr(H)$ be the set of irreducible ordinary
characters and the set of irreducible Brauer characters of $H$,
respectively. If $\chi\in Irr(H)$, it is known that $\chi$ uniquely
defines an algebra homomorphism $\omega_{\chi}: {\bf Z}(\Cc
H)\rightarrow \Cc$ by
$\omega_{\chi}(\hat{C})=\frac{|C|\chi(x_{C})}{\chi(1)}$, where $C$
is a conjugacy class of $H$, $x_{C}\in C$ and $\hat{C}=\sum_{x\in
C}x$. Since $\omega_{\chi}(\hat{C})$ is an algebraic integer, we may
construct an algebra homomorphism $\lambda_{\chi}: {\bf
Z}(FH)\rightarrow F$ by setting
$\lambda_{\chi}(\hat{C})=\omega_{\chi}(\hat{C})^*$. We see that
every irreducible ordinary character $\chi$ gives rise to an algebra
homomorphism ${\bf Z}(FH)\rightarrow F$. This is also true for
irreducible Brauer characters, that is, every irreducible Brauer
character determines an algebra homomorphism ${\bf Z}(FH)\rightarrow
F$. The $2$-{\it blocks} of $H$ are the equivalence classes of
$Irr(H)\cup IBr(H)$ under the equivalence relation $\chi\sim \phi$
if $\lambda_{\chi}=\lambda_{\phi}$ for $\chi$, $\phi\in Irr(H)\cup
IBr(H)$. For basic results on blocks of finite group, we refer the
reader to Chapter 3 of \cite{gabriel}.

The group $H$ has $1$ trivial character $1$ of degree $1$, $2$ 
irreducible ordinary characters $\beta_1$ and $\beta_2$ (respectively, 
$\eta_1$ and $\eta_2)$ of degree $\frac{q+1}{2}$ (respectively, 
$\frac{q-1}{2}$), $1$ irreducible ordinary character $\gamma$ of degree 
$q$, $\frac{q-1}{4}$ (respectively, $\frac{q-3}{4}$) irreducible 
ordinary characters $\chi_s$ for $1\le s \le \frac{q-1}{4}$ 
(respectively, $1\le s \le \frac{q-3}{4}$) of degree $q-1$,
and $\frac{q-5}{4}$ (respectively, $\frac{q-3}{4}$) irreducible
ordinary characters $\phi_r$ of degree $q+1$ for $1\le r\le \frac{q-5}{4}$ 
(respectively, $1\le r \le\frac{q-3}{4}$) if $q\equiv 1\pmod 4$ 
(respectively, $q\equiv 3\pmod 4$). The character table of
$H$ is given in the Appendix. The following lemma describes
how the irreducible ordinary characters of $H$ are partitioned
into $2$-blocks.

\begin{Lemma}\label{blocks}
First assume that $q\equiv 1 \pmod 4$ and $q-1 = m2^n$, where $2\nmid m$.
\begin{itemize}

\item[(i)] The principal block $B_0$ of $H$ contains $2^{n-2}+3$ irreducible characters
$$\chi_0 = 1,\;\gamma,\;\beta_1,\;\beta_2,\;\phi_{i_1},\cdots,\;\phi_{i_{(2^{n-2}-1)}},$$ 
where $\chi_0=1$ is the trivial character of $H$, $\gamma$ is the irreducible character of
degree $q$ of $H$,  $\beta_1$ and $\beta_2$ are the two irreducible characters of degree
$\frac{q+1}{2}$, and $\phi_{i_k}$ for $1\le k\le 2^{n-2}-1$ are distinct irreducible characters
of degree $q+1$ of $H$.

\item[(ii)]  $H$ has $\frac{q-1}{4}$ blocks $B_s$ of defect $0$ for $1\le s \le \frac{q-1}{4}$,
each of which contains an irreducible ordinary character $\chi_s$ of degree $q-1$.

\item[(iii)] If $m\ge 3$, then $H$ has $\frac{m-1}{2}$ blocks $B_t^{'}$ of defect $n-1$ for 
$1\le t\le \frac{m-1}{2}$, each of which contains $2^{n-1}$ irreducible ordinary characters 
$\phi_{t_i}$ for $1\le i \le 2^{n-1}$.

\end{itemize}

Now assume that $q\equiv 3\pmod 4$ and $q+1=m2^n$, where $2\nmid m$ .
\begin{itemize}

\item[(iv)] The principal block $B_0$ of $H$ contains $2^{n-2}+3$ irreducible characters
$$\chi_0 = 1,\;\gamma,\;\eta_1,\;\eta_2,\;\chi_{i_1},\cdots,\;\chi_{i_{(2^{n-2}-1)}},$$
where $\chi_0=1$ is the trivial character of $H$, $\gamma$ is the irreducible character of
degree $q$ of $H$, $\eta_1$ and $\eta_2$ are the two irreducible characters of degree
$\frac{q-1}{2}$, and $\chi_{i_k}$ for $1\le k\le 2^{n-2}-1$ are distinct irreducible characters of
degree $q-1$ of $H$.

\item[(v)] $H$ has $\frac{q-3}{4}$ blocks $B_r$ of defect $0$ for
$1\le r\le \frac{q-3}{4}$, each of which contains an irreducible
ordinary character $\phi_r$ of degree $q+1$.

\item[(vi)] If $m\ge 3$, then $H$ has $\frac{m-1}{2}$ blocks $B_t^{'}$ of 
defect $n-1$ for $1\le t\le \frac{m-1}{2}$, each of which contains $2^{n-1}$
irreducible ordinary characters $\chi_{t_i}$ for $1\le i \le 2^{n-1}$.

\end{itemize}

Moreover, the above blocks form all the $2$-blocks of $H$.
\end{Lemma}

{\Proof} Parts (i) and (iv) are from Theorem 1.3 in {\rm
\cite{landrock}} and their proofs can be found in Chapter 7 of {III}
in {\rm \cite{brauer}}. Parts (ii) and (v) are special cases of
Theorem 3.18 in {\rm \cite{gabriel}}. Parts (iii) and (vi) are
proved in Sections {II} and {VIII} of {\rm \cite{burkhardt}}. \QED


The following result will be used to calculate some of the
block idempotents. Recall that $g\in H$ is $p$-{\it singular}
if $p$ divides the order of $g$.
\begin{Corollary}\label{WBO}\cite[Corollary 3.7]{gabriel}
Suppose that $B$ is a $p$-block of $H$ and let $g,h\in H$.
If $h$ is $p$-regular and $g$ is $p$-singular, then
$$\displaystyle\sum_{\chi \in Irr(B)}\chi(h)\overline{\chi(g)}=0.$$
\end{Corollary}

\subsection{Block Idempotents}

Let $Bl(H)$ be the set of $2$-blocks of $H$. If $B\in Bl(H)$, we write
$$f_B = \displaystyle\sum_{\chi\in Irr(B)}e_{\chi},$$
where $e_{\chi}=\frac{\chi(1)}{|H|}\sum_{g\in H} \chi(g^{-1})g$ is a
central primitive idempotent of $\mathbf{Z}(\Cc H)$ and
$Irr(B)=Irr(H)\cap B$. For future use, we define $IBr(B)=IBr(H)\cap
B$. Since $f_B$ is an element of $\mathbf{Z}(\Cc H)$, we may write
\begin{displaymath}
\begin{array}{lllll}
f_B& =& \displaystyle\sum_{C\in cl(H)}f_B(\widehat{C})\widehat{C},\\
\end{array}
\end{displaymath}
where $cl(H)$ is the set of conjugacy classes of $H$, $\widehat{C}$ is the sum of
elements in the class $C$, and
\begin{equation}\label{id}
\begin{array}{llll}
f_B(\widehat{C})& = &\frac{1}{|H|}\displaystyle\sum_{\chi\in Irr(B)}\chi(1)\chi(x_C^{-1})
\end{array}
\end{equation}
with a fixed element $x_C\in C$.
\begin{Theorem}\label{osima}
Let $B\in Bl(H)$. Then $f_B\in \mathbf{Z}(\mathbf{S}H)$. In other words,
$f_B(\widehat{C})\in \mathbf{S}$ for each block of $H$.

\end{Theorem}
{\Proof} It follows from Corollary 3.8 in {\rm \cite{gabriel}}. \QED

We extend the ring homomorphism $*: \mathbf{S}\rightarrow F$ to a
ring homomorphism $*:\mathbf{S}H\rightarrow FH$ by setting
$(\sum_{g\in H} s_g g)^*= \sum_{g\in H} s_g^* g$. Note that $*$ maps
$\mathbf{Z}(\mathbf{S}H)$ onto $\mathbf{Z}(FH)$ via $(\sum_{C\in
cl(H)}s_C \widehat{C})^*$ = $\sum_{C\in cl(H)} s_C^* \widehat{C}$.
Now we define
$$e_B = (f_B)^* \in \mathbf{Z}(FH),$$
which is the {\it block idempotent} of $B$. Note that $e_B
e_{B^{'}}= \delta_{B B^{'}}e_B$ for $B$, $B^{'}\in Bl(H)$, where
$\delta_{B B^{'}}$ equals 1 if $B=B'$, 0 otherwise. Also
$1=\sum_{B\in Bl(H)}e_B$.

To find $f_B(\widehat{[0]})$ or $e_B(\widehat{[0]})$, we need the following lemma.
\begin{Lemma}\label{p2}
Assume that $q-1=m2^n$ or $q+1=m2^n$ according as $q\equiv 1\pmod 4$
or $q\equiv 3\pmod 4$, where $2\nmid m$. Let $g\in [0]$ and
$\phi_{i_k}$ $($respectively, $\chi_{i_k}$$)$ for $1\le k\le
2^{n-2}-1$ be the ordinary characters of degree $q+1$
$($respectively, $q-1$$)$ in the principal block of $H$ if $q\equiv
1\pmod 4$ $($respectively, $q\equiv 3\pmod 4$$)$ . Then
\begin{displaymath}
\displaystyle\sum_{k=1}^{2^{n-2}-1}\phi_{i_{k}}(g)=
\begin{cases}
-2, &\text{if}\;q\equiv1\pmod 8,\\
0, &\text{if}\;q\equiv 5 \pmod 8,\\
\end{cases}
\end{displaymath}
and
\begin{displaymath}
\displaystyle\sum_{k=1}^{2^{n-2}-1}\chi_{i_{k}}(g)=
\begin{cases}
2, &\text{if}\;q\equiv 7\pmod 8,\\
0, &\text{if}\;q\equiv 3 \pmod 8.\\
\end{cases}
\end{displaymath}
\end{Lemma}
{\Proof} Let $1_H$ be the identity of $H$ and $g\in [0]$. If $q\equiv 1\pmod 4$, 
from Corollary~\ref{WBO}
and the character table of $H$ in the appendix, we have
$$0=\displaystyle\sum_{\chi\in Irr(B_0)}\chi(1_H)\overline{\chi(g)}=(1)(1)+
(q)(1)+ (2)(\frac{q+1}{2})(-1)^{(q-1)/4}+
(q+1)\displaystyle\sum_{k=1}^{2^{n-2}-1}\overline{\phi_{i_{k}}(g)}.$$
Therefore,
\begin{displaymath}
\displaystyle\sum_{k=1}^{2^{n-2}-1}\phi_{i_{k}}(g)=
\begin{cases}
-2, &\text{if}\;q\equiv1\pmod 8,\\
0, &\text{if}\;q\equiv 5 \pmod 8.\\
\end{cases}
\end{displaymath}
The conclusion in the case where $q\equiv 3\pmod 4$ can be proved in the same way.\QED

Using the character tables of $H$ in the Appendix,
Proposition~\ref{classes}, and Lemma~\ref{blocks}, we can find
$e_B(\widehat{C})$ for each $2$-block $B$ and each conjugacy class
$C$ of $H$.
\begin{Lemma}\label{expression}
First assume that $q\equiv 1 \pmod 4$ and $q-1= m2^n$ with $2\nmid m$.
\begin{itemize}

\item[1.] Let $B_0$ be the principal block of $H$. Then
\begin{enumerate}
\renewcommand{\labelenumi}{(\alph{enumi})}
\item 
$e_{B_0}(\widehat{D})=1$.

\item 
 $e_{B_0}(\widehat{F^+})= e_{B_0}(\widehat{F^-})\in F$.

\item 
$e_{B_0}(\widehat{[\theta_i]})\in F$, $e_{B_0}(\widehat{[0]})=0$.

\item 
$e_{B_0}(\widehat{[\pi_k]})=1$.
\end{enumerate}
\item[2.] Let $B_s$ be any block of defect $0$ of $H$. Then
\begin{enumerate}
\renewcommand{\labelenumi}{(\alph{enumi})}
\item 
$e_{B_s}(\widehat{D})=0$.
\item 
$e_{B_s}(\widehat{F^+})= e_{B_s}(\widehat{F^-})=1$.
\item 
$e_{B_s}(\widehat{[0]})= e_{B_s}(\widehat{[\theta_i]})= 0$.
\item 
$e_{B_s}(\widehat{[\pi_k]})\in F$.
\end{enumerate}
\item[3.] Suppose $m\ge 3$ and let $B_t^{'}$ be any block of defect $n-1$ of $H$. Then
\begin{enumerate}
\renewcommand{\labelenumi}{(\alph{enumi})}
\item 
$e_{B_t^{'}}(\widehat{D})=0$.
\item 
$e_{B_t^{'}}(\widehat{F^+})= e_{B_t^{'}}(\widehat{F^-})=1$.
\item 
$e_{B_t^{'}}(\widehat{[\theta_i]})\in F$, $e_{B_t^{'}}(\widehat{[0]})=0$.
\item 
$e_{B_t^{'}}(\widehat{[\pi_k]}) = 0$.
\end{enumerate}
\end{itemize}
Now assume that $q\equiv 3 \pmod 4$. Suppose that $q+1=m2^n$ with $2\nmid m$.
\begin{itemize}

\item[4.] Let $B_0$ be the principal block of $H$. Then
\begin{enumerate}
\renewcommand{\labelenumi}{(\alph{enumi})}
\item 
$e_{B_0}(\widehat{D})=1$.

\item 
$e_{B_0}(\widehat{F^+})= e_{B_0}(\widehat{F^-})\in F$.

\item 
$e_{B_0}(\widehat{[\theta_i]})=1$.

\item 
$e_{B_0}(\widehat{[0]})=0$, $e_{B_0}(\widehat{[\pi_k]})\in F$.
\end{enumerate}
\item[5.] Let $B_r$ be any block of defect $0$ of $H$. Then
\begin{enumerate}
\renewcommand{\labelenumi}{(\alph{enumi})}
\item 
$e_{B_r}(\widehat{D})=0$.
\item 
$e_{B_r}(\widehat{F^+})= e_{B_r}(\widehat{F^-})=1$.
\item 
$e_{B_r}(\widehat{[0]})= e_{B_r}(\widehat{[\pi_k]}) = 0$.
\item 
$e_{B_r}(\widehat{[\theta_i]})\in F$.
\end{enumerate}
\item[6.] Suppose that $m\ge 3$ and let $B_t^{'}$ be any block of defect $n-1$ of $H$. Then
\begin{enumerate}
\renewcommand{\labelenumi}{(\alph{enumi})}
\item 
$e_{B_t^{'}}(\widehat{D})=0$.
\item 
$e_{B_t^{'}}(\widehat{F^+})= e_{B_t^{'}}(\widehat{F^-})=1$.
\item 
$e_{B_t^{'}}(\widehat{[\theta_i]})= 0$.
\item 
$e_{B_t^{'}}(\widehat{[0]})=0$, $e_{B_t^{'}}(\widehat{[\pi_k]})\in F$.
\end{enumerate}

\end{itemize}
\end{Lemma}
{\Proof} We only give the proof for the case where $q\equiv1\pmod 4$.

Since $D$ contains only the identity matrix, by (ii) of
Lemma~\ref{blocks}, (\ref{id}), and the second column of the
character table of $PSL(2,q)$, we have

\begin{equation}
\begin{array}{lllllll}
f_{B_0}(\widehat{D}) &= & \frac{1}{|H|}\displaystyle\sum_{\chi\in Irr(B_0)}
\chi(1)\chi(x_{D}^{-1})\\
{} & = & \frac{1}{qm(\frac{q+1}{2})2^n}[1+q^2+(2)(\frac{q+1}{2})^2 +
(2^{n-2}-1) (q+1)^2]\\
{} & = & \frac{[(2^{2n-2}+2^{n-1})m^2+m2^n+1]}{qm(\frac{q+1}{2})}\\
{}& \in &\Qq
\setminus\Qq\cap\mathcal{P}.\\
\end{array}
\end{equation}
The last inclusion holds since $2\nmid[(2^{2n-2}+2^{n-1})m^2+m2^n+1]$ for $n\ge 2$ and
$2\nmid qm(\frac{q+1}{2})$. Since $f_{B_0}^*(\widehat{D})=e_{B_0}(\widehat{D})$ and $F$
has characteristic $2$, it follows that $e_{B_0}(\widehat{D})=1$. Similarly,
\begin{equation}
\begin{array}{lllllll}
f_{B_0}(\widehat{[\pi_k]}) &= & \frac{1}{|H|}\displaystyle\sum_{\chi\in Irr(B_0)}\chi(1)\chi(x_{[\pi_k]}^{-1})\\
{} & = & \frac{1}{qm(\frac{q+1}{2})2^n}[1+(-1)(q)+(0) (2) (\frac{q+1}{2}) + (0)
(2^{n-2}-1)(q+1)]\\
{} & = & -\frac{m}{qm(\frac{q+1}{2})}\\
{} & \in & \Qq\setminus\Qq\cap\mathcal{P}\\
\end{array}
\end{equation}
and
\begin{equation}\label{p1}
\begin{array}{lllllll}
f_{B_0}(\widehat{[0]}) & = & \frac{1}{|H|}\displaystyle\sum_{\chi\in B_0}\chi(1)\chi(x_{[0]}^{-1})\\
{} & = & \frac{1}{qm(\frac{q+1}{2})2^n}[(q+1)\displaystyle\sum_{k=1}^{2^{n-2}-1}
\phi_{i_k}(x_{[0]}^{-1})+(q+1) (1)+(q+1)(-1)^{(q-1)/4}]\\
{} & = & 0.

\end{array}
\end{equation}
The last equality in (\ref{p1}) follows from Lemma~\ref{p2}.
Therefore, $e_{B_0}([\pi_k])=1$ for $1\le k\le \frac{q-1}{4}$ and
$e_{B_0}([0])=0$. The remaining conclusions in part 1 follow from
Theorem~\ref{osima}.

Let $B_s$ be any block of defect $0$ of $H$. Similar calculations yield
\begin{displaymath}
\begin{array}{lllllll}
f_{B_s}(\widehat{D}) 
\in \mathcal{P},&e_{B_s}(\widehat{D})=0;\\
f_{B_s}(\widehat{F^+})= f_{B_s}(\widehat{F^-}) 
\in\Qq
\setminus\Qq\cap\mathcal{P},& e_{B_s}(\widehat{F^+})=e_{B_s}(\widehat{F^-})=1;\\
f_{B_s}(\widehat{[0]}) =  f_{B_s}(\widehat{[\theta_i]})= 0, & e_{B_s}(\widehat{[0]}) = e_{B_s}(\widehat{[\theta_i]})= 0;\\
f_{B_s}(\widehat{[\pi_k]})\in\mathbf{S}, &e_{B_s}(\widehat{[\pi_k]})\in F. \\
\end{array}
\end{displaymath}

Let $B_t^{'}$ be any block of defect $n-1$ of $H$. Then
\begin{displaymath}
\begin{array}{lllllllll}
f_{B_t^{'}}(\widehat{D})
\in \mathcal{P},&
e_{B_t^{'}}(\widehat{D})=0;\\
f_{B_t^{'}}(\widehat{F^+}) = f_{B_t^{'}}(\widehat{F^-})
\in\Qq\setminus\Qq\cap\mathcal{P},&e_{B_t^{'}}(\widehat{F^+}) = e_{B_t^{'}}(\widehat{F^-})=1; \\
f_{B_t^{'}}(\widehat{[\theta_i]})\in\mathbf{S},\;f_{B_t^{'}}(\widehat{[0]})=0,\;f_{B_t^{'}}
(\widehat{[\pi_k]}) = 0, &e_{B_t^{'}}(\widehat{[\theta_i]})\in F, e_{B_t^{'}}(\widehat{[0]})=0,e_{B_t^{'}}
(\widehat{[\pi_k]}) = 0. \\
\end{array}
\end{displaymath}
\QED

Let $M$ be an $\mathbf{S}H$-module. We denote the reduction $M/\mathcal{P}M$, 
which is an $FH$-module, by $\overline{M}$. Then the following lemma is apparent.
\begin{Lemma}\label{reduction}
Let $M$ be an $\mathbf{S}H$-module and $B\in Bl(H)$. Using the above notation, 
we have $$\overline{f_B M} = e_B \overline{M},$$
i.e. reduction commutes with projection onto a block $B$.
\end{Lemma}

\section{Incidence Matrices and Their Corresponding Maps}

Let $k$ be the complex
field $\Cc$, the algebraic closure $F$ of $\Ff_2$, or the ring 
$\mathbf{S}$ in~(\ref{ring_s}). Let $W$ be a subset of $E$. 
We use ${\bf x}_{W}$ to
denote the characteristic vector of $W$ with respect to $E$; that
is, ${\bf x}_{W}$ is a $(0,1)$-row vector of length $|E|$, whose
entries are indexed by the external points $\p\in E$; the entry of
${\bf x}_W$ indexed by $\p$ is $1$ if and only if $\p\in W$. If
$W=\{\p\}$ is a singleton subset of $E$, then we usually use
${\bf x}_{\p}$ instead of ${\bf x}_{\{\p\}}$ to denote the
characteristic vector of $\{\p\}$ if no confusion occurs.

Let $k^E$ be the free $k$-module with the natural basis
$\{{\bf x}_{\p}\mid \p\in E\}$. It is clear that $k^E$ is a
$kH$-permutation module since $H$ acts on $E$. Let 
$y=\sum_{\p\in E}a_{\p}{\bf x}_{\p}\in k^E$, where $a_{\p}\in k$. 
Then the action of $h\in H$ on $y$ is given by 
$h\cdot y=h\cdot \sum_{\p\in E}
a_{\p}{\bf x}_{\p}=\sum_{\p\in E}a_{\p}(h\cdot {\bf x}_{\p})
=\sum_{\p\in E}a_{\p}{\bf x}_{\p^h}$. Since $H$ is transitive
on $E$, we have
\begin{equation}\label{ind}
k^E\cong \text{Ind}_{K}^H(1_k),
\end{equation}
where $K$ is the stabilizer of an external point in $H$,
$1_k$ is the trivial $kK$-module and $\text{Ind}_K^H(1_k)$ 
is the $kH$-module induced by $1_k$.




Now we define the map
\begin{equation}\label{phi}
\phi: F^E\rightarrow F^E
\end{equation}
by first specifying the images of the natural basis elements under
$\phi$ as follows
$${\bf x}_\p\mapsto \displaystyle\sum_{\q\in \p^\perp\cap E} {\bf x}_\q;$$
then we extend this specification linearly to the whole of $F^E$. As
an $F$-linear map, it is clear that the matrix representation of
$\phi$ with respect to the natural basis (suitably ordered) of $F^E$ 
is $\B$. (Here $\B$ is viewed as a matrix with entries in $F$.) 
Let $\phi^i$ denote the $i$-fold composition of $\phi$. Then it 
is obvious that the matrix of $\phi^i$ with respect to the natural 
basis of $F^E$ is $\B^i$ and $\phi^i({\bf x})={\bf x}\B^i$. Moreover, 
$\phi^5=\phi$ since $\B^5=\B$ by Theorem~\ref{matrix}.

\begin{Lemma}
The above map $\phi$ is an $FH$-module homomorphism from $F^E$ to
$F^E$.
\end{Lemma}
{\Proof} This follows from the fact that $H$ preserves incidence.
\QED
In the rest of the paper, we will always use ${\bf 0}$ and $\hat{{\bf 0}}$
to denote the all-zero row vector of length $|E|$ and the all-zero matrix 
of size $|E|\times |E|$, respectively.
\begin{Proposition}\label{sum}
As $FH$-modules, $F^E=Im(\phi)\oplus Ker(\phi)$, where $Im(\phi)$
and $Ker(\phi)$ are the image and kernel of $\phi$, respectively.
\end{Proposition}
{\Proof} It is clear that $Ker(\phi)\subseteq Ker(\phi^4)$. If
${\bf x}\in Ker(\phi^4)$, then ${\bf x}\in Ker(\phi)$ since
$$\phi({\bf x})=\phi^5({\bf x})=\phi(\phi^4({\bf x}))={\bf 0}.$$
Therefore, $Ker(\phi^4)=Ker(\phi)$. If ${\bf x}\in Ker(\phi)\cap
Im(\phi)$, then there exists a ${\bf y} \in F^E$ such that
${\bf x}=\phi({\bf y})$. It follows that
$${\bf x}=\phi({\bf y})=\phi^5({\bf y})=\phi^4(\phi({\bf y}))=\phi^3(\phi({\bf x}))={\bf 0}.$$
Hence $Ker(\phi)\cap Im(\phi)={\bf 0}$.

It is clear that $Im(\phi)+Ker(\phi)\subseteq F^E$.
Now we assume that ${\bf x}\in F^E$. Note that ${\bf x}=\phi^4({\bf
x})+ {\bf x} -\phi^4({\bf x})$. If ${\bf x}\in Ker(\phi)$, then
${\bf x}\in Im(\phi)+Ker(\phi)$. If ${\bf x}\notin Ker(\phi)$, then
$\phi^4({\bf x})\not={\bf 0}$ since $Ker(\phi^4)=Ker(\phi)$. 
Now $\phi^4({\bf x})\in Im(\phi)$, and 
${\bf x}-\phi^4({\bf x})\in Ker(\phi)$ for
$\phi({\bf x}-\phi^4({\bf x}))=\phi({\bf x})-\phi^5({\bf x})
={\bf 0}$. Therefore ${\bf x}\in Im(\phi)+Ker(\phi)$ and 
thus the proof is complete. \QED

\begin{Corollary}\label{directsum}
As $FH$-modules, $\text{Ind}_K^H(1_F)\cong Ker(\phi)\oplus Im(\phi)$.
\end{Corollary}
{\Proof} The conclusion follows immediately from
Proposition~\ref{sum} and the fact that $\text{Ind}_K^H(1_F)\cong
F^E$. \QED

Next we define
$$\mathbf{C}:=\mathbf{B}^4+I,$$
and
$$\mathbf{D}:=\mathbf{C}+J,$$
where $\mathbf{B}$ is again viewed as a matrix over $F$, $I$ is the
identity matrix and $J$ is the all-one matrix. By
Corollary~\ref{nullspace}, the matrix $\mathbf{C}$ can be viewed as
the incidence matrix between the external points $\p\in E$ and the
subsets ${N_E(\p)}^{\rm a}$ of $E$, $\p\in E$. The matrix $\mathbf{D}$
can also be viewed as an incidence matrix; it is the incidence matrix
between the external points $\p\in E$ and the subsets $N_E^{'}(\p)$ of
$E$, $\p \in E$. Let $\phi_1$ (respectively, $\phi_2$) be the
$FH$-homomorphism from $F^E$ to $F^E$ whose matrix with
respect to the natural basis is $\mathbf{C}$ (respectively, $\mathbf{D}$).

\begin{Lemma}\label{separate}
Assume that $q\equiv 1\pmod 4$. Then as $FH$-modules, we have
$Ker(\phi)=\langle\hat{\mathbf{J}}\rangle\oplus Im(\phi_2)$, where
$\langle\hat{\mathbf{J}}\rangle$ is the trivial $FH$-module
generated by the all one vector $\hat{\bf J}$ of length $|E|$.
\end{Lemma}
{\Proof} Let ${\bf y}\in \langle\hat{\mathbf{J}}\rangle\cap
Im(\phi_2)$. Then ${\bf y}=\phi_2({\bf x})=\lambda\hat{\mathbf{J}}$
for some $\lambda\in F$ and ${\bf x}\in F^E$. By the definition of
$\phi_2$, we have $\phi_2({\bf x})={\bf x}\mathbf{D}$. It follows
that ${\bf x}(\mathbf{B}^4+{I}+{J})=\lambda\hat{\bf J}$ . Note that
${J}^2={J}$ and $\hat{\bf {J}}{J}=\hat{\bf J}$ since $2\nmid |E|$
when $q\equiv 1$ (mod 4). Also each row of $\B^4$ has an
even number of $1$s since the row of $\B^4$ indexed by
$\p\in E$ is the characteristic vector of $N_E(\p)$ or
$N_E(\p)\cup T_E(\p)$ according as $q\equiv 5\pmod 8$ or
$q\equiv 1\pmod 8$ (see the proof of Theorem~\ref{matrix});
that is, $\mathbf{B}^4J={\hat{\bf 0}}$. Thus,
$$\lambda\hat{\bf J}=\lambda\hat{\bf J}J={\bf x}
(\mathbf{B}^4+{I}+{J}){J} = {\bf x}(\mathbf{B}^4{J}+{I}{J}+{J}^2)
={\bf x}(\hat{{\bf 0}}+ {J} + {J}) ={\bf 0}.$$
It follows that $\lambda=0$. We have shown that
$\langle\hat{\mathbf{J}}\rangle\cap Im(\phi_2)={\bf 0}$.

It is obvious that $\langle\hat{\mathbf{J}}\rangle +
Im(\phi_2)\subseteq Ker(\phi)$. Let ${\bf x}\in Ker(\phi)$. Then by
Corollary~\ref{nullspace}, there exists a ${\bf y}\in F^E$ such that
${\bf x}={\bf y}(\mathbf{B}^4+{I})$, which in turn is equal to ${\bf
y}(\mathbf{B}^4+{I}+ {J})+\langle{\bf y}, \hat{{\bf
J}}\rangle\hat{{\bf J}}$ since ${\bf y}{J}=\langle{\bf y}, \hat{{\bf
J}}\rangle\hat{{\bf J}}$, where $\langle{\bf y}, \hat{{\bf
J}}\rangle$ is the standard inner product of the vectors ${\bf y}$
and $\hat{{\bf J}}$. Hence ${\bf x}\in
\langle\hat{\mathbf{J}}\rangle + Im(\phi_2)$ and thus
$Ker(\phi)=\langle\hat{\mathbf{J}}\rangle\oplus Im(\phi_2)$. \QED

\section{An Induced Character}
In this section, we consider the induced complex character
$1\uparrow_{K}^H$ afforded by the $\Cc H$-module
$\text{Ind}_K^H(1_{\Cc})$ and decompose $1\uparrow_K^H$
into the sum of irreducible complex characters of $H$ by
using the well-known Frobenius reciprocity \cite{frob}, the 
character tables of $H$, and Lemma~\ref{stabilizer_intersection}.

\begin{Lemma}\label{decomposition_1}
Assume that $q\equiv 1\pmod 4$. Let $\chi_s$,
$1\le s \le \frac{q-1}{4}$, be the irreducible ordinary
characters of degree $q-1$, $\phi_r$, $1\le r \le \frac{q-5}{4}$,
irreducible ordinary characters of degree $q+1$, $\gamma$
the irreducible of degree $q$, and $\beta_j$, $1\le j\le 2$,
irreducible ordinary characters of degree $\frac{q+1}{2}$.

\begin{itemize}
\item[(i)] If $q\equiv 1 \pmod 8$, then $$1\uparrow_{K}^H =
1 + \displaystyle\sum_{s=1}^{(q-1)/4}\chi_s + 2\gamma +
\beta_1+\beta_2 + \displaystyle\sum_{j=1}^{(q-9)/4}\phi_{r_j},$$
where $\phi_{r_j}$, $1\le j\le\frac{q-9}{4}$, may not be distinct.
\item[(ii)] If $q\equiv 5\pmod 8$, then $$1\uparrow_{K}^H =
1 + \displaystyle\sum_{s=1}^{(q-1)/4}\chi_s + 2\gamma +
\displaystyle\sum_{j=1}^{(q-5)/4}\phi_{r_j},$$ where
$\phi_{r_j}$, $1\le j\le\frac{q-5}{4}$, may not be distinct.
\end{itemize}

Next assume that $q\equiv 3\pmod 4$. Let $\chi_s$, $1\le s \le
\frac{q-3}{4}$, be the irreducible ordinary characters of degree
$q-1$, $\phi_r$, $1\le r \le \frac{q-3}{4}$, the irreducible ordinary
characters of degree $q+1$, $\gamma$ the irreducible
character of degree $q$, and $\eta_j$, $1\le j\le 2$, the
irreducible ordinary characters of degree $\frac{q-1}{2}$.

\begin{itemize}
\item[(iii)] If $q\equiv 3\pmod 8$, then $$1\uparrow_{K}^H=
1+\displaystyle\sum_{r=1}^{(q-3)/4}\phi_r+\gamma+\eta_1+
\eta_2+\displaystyle\sum_{j=1}^{(q-3)/4}\chi_{s_j},$$
where $\chi_{s_j}$, $1\le j\le\frac{q-3}{4}$, may not be distinct.

\item[(iv)] If $q\equiv 7\pmod 8$, then $$1\uparrow_{K}^H=1
+\displaystyle\sum_{r=1}^{(q-3)/4}\phi_r+\gamma+
\displaystyle\sum_{j=1}^{(q+1)/4}\chi_{s_j},$$ where
$\chi_{s_j}$, $1\le j\le\frac{q+1}{4}$, may not be distinct.
\end{itemize}
\end{Lemma}
{\Proof} We only give the detailed proof for the case where
$q\equiv 1\pmod 4$.

Let $1_H$ be the trivial character of $H$. By the Frobenius
reciprocity \cite{frob},
$$\left\langle 1\uparrow_{K}^H, 1_H\right\rangle_H=
\left\langle1, 1_H\downarrow_K^H\right\rangle_K = 1.$$

Let $\chi_s$ be an irreducible character of degree $q-1$
of $H$, where $1\le s \le \frac{q-1}{4}$. We denote the
number of elements of $K$ lying in the class $[\pi_k]$ by
$d_k$. Then $d_k=0$
by Lemma~\ref{stabilizer_intersection}.
\begin{displaymath}
\begin{array}{lllll}
\left\langle1\uparrow_K^H,\chi_s\right\rangle_H & =
& \left\langle1, \chi_s\downarrow_K^H\right\rangle_K\\
{}&= & \frac{1}{|K|}\displaystyle\sum_{g\in K}\chi_s
\downarrow_{K}^H(g)\\
{} & = & \frac{1}{q-1}[(1)(q-1)+\displaystyle\sum_{k=0}^{(q-1)/4}
(-d_k\delta^{(2k)s}-d_k\delta^{-(2k)s})]\\
{}& = & 1\\
\end{array}
\end{displaymath}

Let $\gamma$ be the irreducible character of degree $q$ of $H$.
\begin{displaymath}
\begin{array}{llllll}
\left\langle1\uparrow_K^H,\gamma \right\rangle_H & = &
\left\langle1, \gamma\downarrow_K^H\right\rangle_K \\
{} & = & \frac{1}{|K|}\displaystyle\sum_{g\in K} \gamma
\downarrow_K^H(g)\\
{}& = & \frac{1}{q-1}\left[(1)(q)+(1)(2)(\frac{q-5}{4})+
(1)(\frac{q+1}{2})\right]\\
{} & = & 2\\
\end{array}
\end{displaymath}

Let $\beta_j$ be any irreducible character of degree
$\frac{q+1}{2}$ of $H$.
\begin{equation}\label{charsum1}
\begin{array}{llll}
\left\langle1\uparrow_K^H, \beta_j\right\rangle_H & =
& \frac{1}{|K|}\displaystyle\sum_{g\in K}\beta_j\downarrow_K^H(g)\\
{} & = & \frac{1}{q-1}[(1)(\frac{q+1}{2})+(2)\displaystyle
\sum_{i=1}^{(q-5)/4}\zeta(\theta_i)+(\frac{q+1}{2})(-1)^{(q-1)/4}].\\
\end{array}
\end{equation}
If $q\equiv 1\pmod 8$, then
$$\frac{q+1}{2}+(2)\displaystyle\sum_{i=1}^{(q-5)/4}\zeta(\theta_i)+
(\frac{q+1}{2})(-1)^{(q-1)/4}=(q+1)+(2)\displaystyle
\sum_{i=1}^{(q-5)/4}\zeta(\theta_i).$$
By the fact that $\zeta(\theta_i)=1$ or $-1$, we have
$$\frac{q+7}{2}\le (q+1)+(2)\displaystyle\sum_{i=1}^{(q-5)/4}
\zeta(\theta_i)\le \frac{3(q-1)}{2}.$$ Since
$\left\langle1\uparrow_K^H, \beta_j\right\rangle_H$ is a
non-negative integer and $q-1$ is the only integer in the interval
$[\frac{q+7}{2},\frac{3(q-1)}{2}]$ that is divisible by $q-1$ when
$q\ge 9$, we must have that
$$\displaystyle\sum_{i=1}^{(q-5)/4}\zeta(\theta_i)=-1.$$ Similarly,
if $q\equiv 5\pmod 8$, then
$$\displaystyle\sum_{i=1}^{(q-5)/4}\zeta(\theta_i)=0.$$ Therefore,
\begin{displaymath}
\left\langle1\uparrow_K^H, \beta_j\right\rangle_H=
\begin{cases}
1, & \text{if}\; q\equiv 1\pmod 8,\\
0, & \text{if}\; q\equiv 5\pmod 8.\\
\end{cases}\\
\end{displaymath}

Since the sum of the degrees of $1$, $\chi_s$, $\gamma$, and
$\beta_j$ is less than the degree of $1\uparrow_K^H$ and only
the irreducible characters of degree $q+1$ of $H$ have not been
taken into account yet, we see that all the irreducible constituents of
$$1\uparrow_K^H - 1 - \displaystyle\sum_{s=1}^{(q-1)/4}
\chi_s-2\gamma-\beta_1-\beta_2$$
or
$$1\uparrow_K^H-1-\displaystyle\sum_{s=1}^{(q-1)/4}
\chi_s-2\gamma$$
must have degree $q+1$. \QED
\begin{Corollary}\label{char11}
Using the above notation,  
\begin{itemize}
\item[(i)] if $q\equiv 1\pmod 4$, then the character of $f_{B_s}\cdot Ind_K^H (1_{\Cc})$ 
is $\chi_s$ for each block $B_s$ of defect $0$;
\item[(ii)] if $q\equiv 3\pmod 4$, then the character of $f_{B_r}\cdot Ind_K^H (1_{\Cc})$
is $\phi_r$ for each block $B_r$ of defect $0$. 
\end{itemize}
\end{Corollary}
{\Proof} The corollary follows from Lemma~\ref{blocks} and Lemma~\ref{decomposition_1}.
\QED
\section{Statement and Proof of Main Theorem}
We state and prove the main theorem in this section.

\begin{Lemma}\label{y4}
Let $e_{B}$ for $B\in Bl(H)$ be all the primitive
idempotents of $\mathbf{Z}(FH)$. Assume $q-1=2^n m$ or $q+1=2^n m$
with $2\nmid m$ depending on whether $q\equiv 1\pmod 4$ 
or $q\equiv 3\pmod 4$. Using the above notation, 
\begin{itemize}
\item[(i)] if $q\equiv 1\pmod 4$, then $e_{B_0}Im(\phi_2)={\bf 0}$,
$e_{B_s}Im(\phi) ={\bf 0}$ for $1\le s\le \frac{q-1}{4}$, and $e_{{B_t}^{'}}Im(\phi_2)={\bf 0}$
for $m\ge 3$ and $1\le t \le \frac{m-1}{2}$;
\item[(ii)] if $q\equiv 3\pmod 4$, then $e_{B_0}Ker(\phi)={\bf 0}$, 
$e_{B_r}Im(\phi) ={\bf 0}$ for $1 \le r \le \frac{q-3}{4}$, and $e_{B_t^{'}}Ker(\phi) ={\bf 0}$ 
for $m\ge 3$ and $1\le t\le \frac{m-1}{2}$.
\end{itemize}

\end{Lemma}

{\Proof} It is clear that $\{\x_\p\B\mid \p\in E\}$, $\{\x_\p\C\mid \p\in E\}$,
and $\{\x_\p\D\mid \p\in E\}$ span $Im(\phi)$, $Ker(\phi)$, and
$Im(\phi_2)$ over $F$, respectively. Also, ${\x_\p}\B$, ${\x_\p}\C$, and
${\x_\p}\D$ are the characteristic vectors of $E_{\p^\perp}$, $E_{N_E(\p)^{\rm a}}$,
$E_{N_E(\p)^{'}}$, respectively. Let $B\in Bl(H)$. We notice that
\begin{equation*}
\begin{array}{llllll}
e_{B}{\x_{\p^\perp}} &= &\displaystyle\sum_{C\in cl(H)}e_{B}
(\widehat{C})\displaystyle\sum_{h\in C}h\cdot {\x}_{\p^\perp}\\
{} & = & \displaystyle\sum_{C\in cl(H)}e_{B}(\widehat{C})
\displaystyle\sum_{h\in C}{\x}_{(\p^\perp)^h},\\
{} & = & \displaystyle\sum_{C\in cl(H)}e_{B}(\widehat{C})
\displaystyle\sum_{h\in C}\sum_{\q\in(\p^{\perp})^h\cap E}{\bf x}_{\q};
\end{array}
\end{equation*}
that is,
\begin{equation*}
e_{B}{\x_{\p^\perp}} =\sum_{\q\in E}\mathcal{S}_1(B,\p,\q){\bf x}_{\q},
\end{equation*}
where $\mathcal{S}_1(B,\p,\q):=\sum_{C\in cl(H)}
|C\cap\mathcal{H}_{\p,\q}|e_{B}(\widehat{C})$.
Similarly, we have that
\begin{equation*}
e_{B}{\x_{N_{E}(\p)^{\rm a}}} =\sum_{\q\in E}\mathcal{S}_2(B,\p,\q){\bf x}_{\q}
\;\text{and}\;
e_{B}{\x_{N_{E}(\p)^{'}}} =\sum_{\q\in E}\mathcal{S}_3(B,\p,\q){\bf x}_{\q},
\end{equation*}
where $$\mathcal{S}_2(B,\p,\q):=\sum_{C\in cl(H)}
|C\cap\mathcal{U}_{\p,N_E(\q)^{\rm a}}|e_{B}(\widehat{C})$$
and
$$\mathcal{S}_3(B,\p,\q):=\sum_{C\in cl(H)}
|C\cap\mathcal{U}_{\p,N_E(\q)^{'}}|e_{B}(\widehat{C}).$$

First we assume that $q\equiv 1\pmod 4$.

If $\q=\p$, since
$|C\cap\mathcal{H}_{\p,\p}|$ and $|C\cap\mathcal{U}_{\p,N_E(\p)^{'}}|$ 
is always even for $C\not=[0]$ by Lemma~\ref{pisq} and Lemma~\ref{y3} 
(i) respectively, and $e_{B_0}(\widehat{[0]})=e_{B_s}(\widehat{[0]})=e_{B_t^{'}}(\widehat{[0]})=0$ by 1(c), 2(c), 3(c) of
Lemma~\ref{expression}, then $\mathcal{S}_1(B_s,\p,\p)=
\mathcal{S}_3(B_0,\p,\p)=\mathcal{S}_3(B_t^{'},\p,\p)=0$. 
If $\ell_{\p,\q}\in Pa_{\p}$, then $\mathcal{S}_1(B_s,\p,\q)=0$
since $|C\cap\h_{\p,\q}|$ is always even by 
Lemma~\ref{line_1} (iii) for $C\not=[0]$ and $e_{B_s}(\widehat{[0]})=0$;
by Lemma~\ref{y_1} (i), 1(a), (d) and 3(a), (d) of Lemma~\ref{expression},
$$\mathcal{S}_3(B_0,\p,\q)=e_{B_0}(\widehat{D})+e_{B_0}(\widehat{[\pi_k]})
=1+1=0$$
and
$$\mathcal{S}_3(B_t^{'},\p,\q)=e_{B_t^{'}}(\widehat{D})+e_{B_t^{'}}(\widehat{[\pi_k]})
=0+0=0.$$
If $\q\not=\p$, $\ell_{\p,\q}\in Se_{\p}$,
and $\q\notin \p^\perp$, by Lemma~\ref{line_1} (i) and 2(c) of
Lemma~\ref{expression}, then
$$\mathcal{S}_1(B_s,\p,\q)=e_{B_s}(\widehat{[\theta_{i_1}]})+e_{B_s}
(\widehat{[\theta_{i_2}]})+|[0]\cap\h_{\p,\q}|e_{B_s}(\widehat{[0]})=0+0+0=0;$$
by Lemma~\ref{y_1} (ii) and 1(c), 3(c) of Lemma~\ref{expression},
$\mathcal{S}_3(B_0,\p,\q)=\mathcal{S}_3(B_t^{'},\p,\q)=0$.
If $\q\in\p^\perp$ and $\ell_{\p,\q}\in Se_{\p}$, by
Lemma~\ref{line_1} (ii) and 2(a), (c) of Lemma~\ref{expression},
$$\mathcal{S}_1(B_s,\p,\q)=e_{B_s}(\widehat{D})+|[0]\cap\h_{\p,\q}
|e_{B_s}(\widehat{[0]})=0+0=0;$$ by Lemma~\ref{y_1} (ii) and 1
(c), 3(c) of Lemma~\ref{expression},
$\mathcal{S}_3(B_0,\p,\q)=\mathcal{S}_3(B_t^{'},\p,\q)=0$.
In the case where $\ell_{\p,\q}\in T_{\p}$ and $\q\not=\p$, by
Lemma~\ref{line_1} (iv) and 2(c) of Lemma~\ref{expression}, 
$$\mathcal{S}_1(B_s,\p,\q)=e_{B_s}(\widehat{[\theta_{i_1}]})
+\cdots+e_{B_s}(\widehat{[\theta_{i_m}]})+|[0]\cap\h_{\p,\q}
|e_{B_s}(\widehat{[0]})=0+\cdots+0+0=0$$ for some $m\ge 1$;
if $q\equiv 1\pmod 8$,
by Lemma~\ref{y_1} (iii) and 1(c), (d) and 3 (c), (d),
$$\mathcal{S}_3(B_0,\p,\q)=\displaystyle\sum_{k=1}^{(q-1)/4}e_{B_0}
(\widehat{[\pi_k]})+|[0]\cap\h_{\p,\q}|e_{B_0}(\widehat{[0]})=
\displaystyle\sum_{k=1}^{(q-1)/4}1+0=0,$$
$$\mathcal{S}_3(B_t^{'},\p,\q)=\displaystyle\sum_{k=1}^{(q-1)/4}e_{B_t^{'}}(\widehat{[\pi_k]})+
|[0]\cap\h_{\p,\q}|e_{B_t^{'}}(\widehat{[0]})=\displaystyle\sum_{k=1}^{(q-1)/4}0+0=0;$$
if $q\equiv 5\pmod 8$
by Lemma~\ref{y_1} (iii) and 1(a), (c), (d) and 3 (a), (c), (d),
$$\mathcal{S}_3(B_0,\p,\q)=\displaystyle\sum_{k=1}^{(q-1)/4}e_{B_0}
(\widehat{[\pi_k]})+e_{B_0}(\widehat{D})+|[0]\cap\h_{\p,\q}|e_{B_0}(\widehat{[0]})
=\displaystyle\sum_{k=1}^{(q-1)/4}1+1+0=0,$$
$$\mathcal{S}_3(B_t^{'},\p,\q)=\displaystyle\sum_{k=1}^{(q-1)/4}e_{B_t^{'}}(\widehat{[\pi_k]})+
+e_{B_t^{'}}(\widehat{D})+|[0]\cap\h_{\p,\q}|e_{B_t^{'}}(\widehat{[0]})=
\displaystyle\sum_{k=1}^{(q-1)/4}0+0+0=0.$$
So we have shown that $\mathcal{S}_1(B_s, \p,\q)=0$, $\mathcal{S}_3(B_0,\p,\q)=0$,
and $\mathcal{S}_3(B_t^{'},\p,\q)=0$ for $\p,\q\in E$. The proof of part (i) is completed.

Now we assume that $q\equiv 3\pmod 4$.

If $\q=\p$, since
$|C\cap\mathcal{H}_{\p,\p}|$ and $|C\cap\mathcal{U}_{\p,N_E(\p)^{\rm a}}|$ 
is always even for $C\not=[0]$ by Lemma~\ref{pisq} and Lemma~\ref{y3} 
(i) respectively, and $e_{B_r}(\widehat{[0]})=e_{B_0}(\widehat{[0]})=
e_{B_t^{'}}(\widehat{[0]})=0$ by 4(c), 5(c), 6(d) of
Lemma~\ref{expression}, then $\mathcal{S}_1(B_r,\p,\p)=
\mathcal{S}_2(B_0,\p,\p)=\mathcal{S}_2(B_t^{'},\p,\p)=0$. 
If $\ell_{\p,\q}\in Pa_{\p}$ and $\q\notin \p^\perp$, by Lemma~\ref{y_22} (iii) and
5(c) of Lemma~\ref{expression}, we have
$$\mathcal{S}_1(B_r,\p,\q)=e_{B_r}(\widehat{[\pi_{k_1}]})+e_{B_r}(\widehat{[\pi_{k_2}]})+
|[0]\cap\h_{\p,\q}|e_{B_r}(\widehat{[0]})=0+0+0=0;$$
by Lemma~\ref{y_2} (ii) and
4(d), 6(d) of Lemma~\ref{expression}, we have
$$\mathcal{S}_2(B_0,\p,\q)=|[0]\cap\mathcal{U}_{\p,N_E(\q)^{\rm a}}|e_{B_0}(\widehat{[0]})=0$$ and
$$\mathcal{S}_2(B_t^{'},\p,\q)=|[0]\cap\mathcal{U}_{\p,N_E(\q)^{\rm a}}|e_{B_t^{'}}(\widehat{[0]})=0.$$
If $\ell_{\p,\q}\in Pa_{\p}$ and $\q\in \p^\perp$, by Lemma~\ref{y_22} (ii)
and 2(a) of Lemma~\ref{expression}, we have
$\mathcal{S}_1(B_r,\p,\q)=e_{B_r}(\widehat{[0]})+|[0]\cap\h_{\p,\q}|e_{B_r}(\widehat{[0]})=0$;
by Lemma~\ref{y_2} (ii) and 4(d), 6(d) of Lemma~\ref{expression}, 
$\mathcal{S}_2(B_0,\p,\q)=\mathcal{S}_2(B_t^{'},\p,\q)=0$. 
In the case where $\ell_{\p,\q}\in T_{\p}$, by
Lemma~\ref{y_22} (iv) and 5(c) of Lemma~\ref{expression},
we always have $\mathcal{S}_1(B_r, \p,\q)=0$; if $q\equiv 3\pmod 8$,
by Lemma~\ref{y_2} (iii) and 4(a), 4(c) and 6(a), 4(c) of Lemma~\ref{expression},
$$\mathcal{S}_2(B_0, \p,\q)=\displaystyle\sum_{i=1}^{(q-3)/4}e_{B_0}(\widehat{[\theta_i]})
+|[0]\cap\mathcal{U}_{\p,N_E(\q)^{\rm a}}| e_{B_0}(\widehat{[0]})=1+\displaystyle\sum_{i=1}^{(q-3)/4}1+0=0$$
and
$$\mathcal{S}_2(B_t^{'}, \p,\q)=\displaystyle\sum_{i=1}^{(q-3)/4}e_{B_t^{'}}(\widehat{[\theta_i]})
+|[0]\cap\mathcal{U}_{\p,N_E(\q)^{\rm a}}| e_{B_t^{'}}(\widehat{[0]})=\displaystyle\sum_{i=1}^{(q-3)/4}0+0=0;$$
if $q\equiv 7\pmod 8$,
by Lemma~\ref{y_2} (iii) and 4(c), 6(c) of Lemma~\ref{expression},
$$\mathcal{S}_2(B_0, \p,\q)=e_{B_0}(\widehat{D})+\displaystyle\sum_{i=1}^{(q-3)/4}e_{B_0}(\widehat{[\theta_i]})
+|[0]\cap\mathcal{U}_{\p,N_E(\q)^{\rm a}}| e_{B_0}(\widehat{[0]})=1+\displaystyle\sum_{i=1}^{(q-3)/4}1+0=0$$
and
$$\mathcal{S}_2(B_t^{'}, \p,\q)=e_{B_t^{'}}(\widehat{D})+\displaystyle\sum_{i=1}^{(q-3)/4}e_{B_t^{'}}(\widehat{[\theta_i]})
+|[0]\cap\mathcal{U}_{\p,N_E(\q)^{\rm a}}| e_{B_t^{'}}(\widehat{[0]})=0+\displaystyle\sum_{i=1}^{(q-3)/4}0+0=0.$$
So we have shown that $\mathcal{S}_1(B_r, \p,\q)=0$, $\mathcal{S}_2(B_0,\p,\q)=0$,
and $\mathcal{S}_2(B_t^{'},\p,\q)=0$ for $\p,\q\in E$. The proof of part (ii) is finished.

\QED

\begin{Theorem}\label{maintheorem}
Let $Ker(\phi)$ be the kernel of $\phi$ defined as above. 
\begin{itemize}
\item[(i)] If $q\equiv1\pmod 4$, then
$$Ker(\phi)=\langle{\bf \hat{J}}\rangle\oplus(\displaystyle\bigoplus_{s=1}^{(q-1)/4}N_s),$$
where $\langle{\bf \hat{J}}\rangle$ is the trivial $FH$-module and
$N_s$ for $1\le s\le \frac{q-1}{4}$ are pairwise nonisomorphic projective
simple $FH$-modules of dimension $q-1$; 
\item[(ii)] if $q\equiv 3 \pmod 4$, then
$$Ker(\phi)= \displaystyle\bigoplus_{r=1}^{(q-3)/4}N_r,$$
where $N_r$ for $1\le r\le \frac{q-3}{4}$ are pairwise nonisomorphic
projective simple $FH$-modules of dimension $q+1$.
\end{itemize}
\end{Theorem}
{\Proof}  We have $F^E=Ker(\phi)\oplus Im(\phi)$ by Proposition~\ref{sum}.
Let $B$ be a $2$-block of defect zero of $H$.
Since $e_{B}Im(\phi)={\bf 0}$, by Lemma~\ref{reduction} we have
$$e_{B}Ker(\phi)=e_{B}F^E= \overline{f_{B}{\bf S}^E}.$$
Therefore by Corollary~\ref{char11}, $e_{B}Ker(\phi)=N$,
where $N$ is the projective simple module in $B$.

Suppose $q\equiv 3 \pmod 4$. By Lemma~\ref{y4} (ii), we have
$e_{B_0}Ker(\phi)={\bf 0}$ and $e_{B'_t}Ker(\phi)={\bf 0}$ for all $t$.
Thus, $$Ker(\phi)=\displaystyle\bigoplus_{r=1}^{(q-3)/4} e_{B_r}Ker(\phi)
=\displaystyle\bigoplus_{r=1}^{(q-3)/4}N_r,$$
where $N_r$ is the projective simple module in $B_r$ by the discussion
in the first paragraph.
The theorem is proved in the case where $q\equiv 3 \pmod 4$.

Suppose $q\equiv 1 \pmod 4$.
By Lemma~\ref{y4} (ii), since $e_{B_0}Im(\phi_2)={\bf 0}$
and $e_{B'_t}Im(\phi_2)={\bf 0}$ for all $t$,  we have
$$Im(\phi_2)=\displaystyle\bigoplus_{s=1}^{(q-1)/4} e_{B_s}Im(\phi_2).$$
By Lemma~\ref{separate}, $Ker(\phi)=\langle{\hat{\bf J}}\rangle\oplus Im(\phi_2)$, so since 
$e_{B_s}\langle{\hat{\bf J}}\rangle={\bf 0}$, we have from the first paragraph,
$$
e_{B_s}Im(\phi_2)=e_{B_s}Ker(\phi)=e_{B_s}F^E=N_s,
$$
where $N_s$ is the simple module in $B_s$.
Thus, by Lemma~\ref{separate}, we have
$$Ker(\phi)=\langle{\hat{\bf J}}\rangle\oplus(\bigoplus_{s=1}^{(q-1)/4}N_s),$$
and the theorem is proved.
\QED

The following corollary is immediate.
\begin{Corollary}
The dimension of the code $\mathcal{L}$ generated by the null space of $\B$ over $\Ff_2$ is
\begin{displaymath}
dim_{\Ff_2}(\mathcal{L})=
\begin{cases}
\frac{(q-1)^2}{4} + 1, & \text{if}\;q\equiv 1\pmod 4,\\
\frac{(q-1)^2}{4}  - 1, & \text{if}\;q\equiv 3\pmod 4.\\
\end{cases}
\end{displaymath}
\end{Corollary}
\newpage
\section*{APPENDIX}
The character tables of $PSL(2,q)$ were obtained by Jordan \cite{jordan}
and Schur \cite{schur} independently, from which we can
deduce the character tables of $H$ as follows.
Let $\epsilon\in \Cc$ be a primitive $(q-1)$-th root of
unity and $\delta\in \Cc$ a primitive $(q+1)$-th root of 
unity.

\begin{table}[htp]
\caption{Character table of $H$ when $q\equiv 1\pmod 4$}
\bigskip
\begin{tabular}{|c|c|c|c|c|c|c|c|}
\hline
 Number & 1 & 2 & $\frac{q-5}{4}$ & 1 & $\frac{q-1}{4}$ \\
\hline
 Size & 1 & $\frac{q^2-1}{2}$ & $q(q+1)$ & $\frac{q(q+1)}{2}$ & $q(q-1)$ \\
\hline
Representative &  $D$ & $F^{\pm}$ & $[\theta_i]$ & $[0]$ & $[\pi_k]$  \\
\hline
$\phi_{r} $ & $q+1$ & 1 & $\epsilon^{(2i)r}+\epsilon^{-(2i)r}$ & $2(-1)^{r}$& 0 \\

$\gamma$ & $q$ & 0 & 1& 1 & $-1$ \\

$1$ & 1 & 1 & 1 & 1 & 1 \\

$\chi_s$ & $q-1$ & $-1$ & 0 & 0 & $-\delta^{(2k)s}-\delta^{-(2k)s}$ \\

$\beta_1$ & $\frac{q+1}{2}$ & $\frac{1}{2}(1\pm\sqrt{q})$ & $\zeta(\theta_i)$ & $(-1)^{(q-1)/4}$ & 0 \\

$\beta_2$ & $\frac{q+1}{2}$ & $\frac{1}{2}(1\mp \sqrt{q})$ & $\zeta(\theta_i)$ & $(-1)^{(q-1)/4}$ & 0 \\

\hline
\end{tabular}\label{table:tab6.2}
\end{table}

Here $s=1,2,...,\frac{q-1}{4}$, $r=1,2,...,\frac{q-5}{4}$, $k=1,2,...,\frac{q-1}{4}$, $i=1,2,...,\frac{q-5}{4}$, and $\zeta(\theta_i) =1$ or $-1$. 

\begin{table}[htp]
\caption{Character table of $H$ when $q\equiv 3\pmod 4$}
\bigskip
\begin{tabular}{|c|c|c|c|c|c|c|c|}
\hline
 Number & 1 & 2 & $\frac{q-3}{4}$ & 1 & $\frac{q-3}{4}$ \\
\hline
 Size & 1 & $\frac{q^2-1}{2}$ & $q(q+1)$ & $\frac{q(q-1)}{2}$ & $q(q-1)$ \\
\hline
Representative &  $D$ & $F^{\pm}$ & $[\theta_i]$ & $[0]$ & $[\pi_k]$  \\
\hline
$\phi_r$ & $q+1$ & 1 & $\epsilon^{(2i)r}+\epsilon^{-(2i)r}$ & 0 & 0 \\

$\gamma$ & $q$ & 0 & 1 & $-1$ & $-1$ \\

$1$ & 1 & 1 & 1 & 1 & 1 \\

$\chi_s$ & $q-1$ & $-1$ & 0 & $-2(-1)^{s}$ &$-\delta^{(2k)s}-\delta^{-(2k)s}$ \\

$\eta_1$ & $\frac{q-1}{2}$ & $\frac{1}{2}(-1\pm\sqrt{-q})$ & $0$ & $(-1)^{(q+5)/4}$
&$-\zeta(\pi_k)$ \\

$\eta_2$ & $\frac{q-1}{2}$ & $\frac{1}{2}(-1\mp\sqrt{-q})$ & $0$ & $(-1)^{(q+5)/4}$
&$-\zeta(\pi_k)$ \\
\hline
\end{tabular}\label{table:tab6.4}
\end{table}

Here $s=1,2,...,\frac{q-3}{4}$, $r=1,2,...,\frac{q-3}{4}$, $k=1,2,...,\frac{q-3}{4}$, $i=1,2,...,\frac{q-3}{4}$, and $\zeta(\pi_k) =1$ or $-1$.

\end{document}